\RequirePackage{fix-cm}
\documentclass[smallextended,numbook]{svjour3}       % onecolumn (second format)
\smartqed  % flush right qed marks, e.g. at end of proof
\usepackage{amssymb}
\usepackage{amsmath}
\usepackage{amsfonts}

\DeclareMathOperator{\sgn}{sgn}

\usepackage{algorithm,algorithmic}
\floatname{Algorithm}{Procedure}
\renewcommand{\algorithmicrequire}{\textbf{Input:}}
\renewcommand{\algorithmicensure}{\textbf{Output:}}
\usepackage{enumerate}
\usepackage{multirow}
\usepackage{graphicx}
\graphicspath{{figures/}}
\usepackage{subfigure}
\usepackage{float}
\usepackage{tikz}
\usetikzlibrary{arrows,decorations.pathmorphing,backgrounds,positioning,fit,petri,graphs}
\usepackage{epstopdf}
\ifpdf
  \DeclareGraphicsExtensions{.eps,.pdf,.png,.jpg}
\else
  \DeclareGraphicsExtensions{.eps}
\fi
\usepackage{cite}
\usepackage[colorlinks, linkcolor=blue, anchorcolor=green, citecolor=black]{hyperref}
\usepackage{cleveref}
\usepackage{enumerate}
\usepackage[shortlabels]{enumitem}
\usepackage{orcidlink}
\usepackage{bm}
\DeclareMathAlphabet{\mathpzc}{OT1}{pzc}{m}{it}
\newcommand{\setsep}{\,|\,}
\newcommand{\sph}{\mathbb{S}}
\newcommand{\R}{\mathbb{R}}
\newcommand{\C}{\mathbb{C}}

\newcommand{\Z}{\mathbb{Z}}

\begin{document}

\title{Spherical Designs for Function Approximation and Beyond %\thanks{Grants or other notes
%about the article that should go on the front page should be
%placed here. General acknowledgments should be placed at the end of the article.}
}
%\subtitle{Do you have a subtitle?\\ If so, write it here}

\titlerunning{Spherical Designs for Function Approximation and Beyond}        % if too long for running head

\author{Yuchen Xiao \orcidlink{0000-0003-3844-0952}   \and Xiaosheng Zhuang \orcidlink{0000-0001-7238-0143}
         %etc.
}

%\authorrunning{Short form of author list} % if too long for running head

\institute{Y. Xiao and X. Zhuang (Corresponding Author) \at
              Department of Mathematics, City University of Hong Kong, Tat Chee Avenue, Kowloon Tong, Hong Kong SAR, China.
              \email{yc.xiao@my.cityu.edu.hk; xzhuang7@cityu.edu.hk}           %  \\
              }

\date{Received: date / Accepted: date}
% The correct dates will be entered by the editor

\maketitle

% REQUIRED
\begin{abstract}
In this paper, we compare two optimization algorithms using full Hessian and approximation Hessian to obtain numerical spherical designs through their variational characterization. Based on the obtained spherical design point sets, we investigate the approximation of smooth and non-smooth functions by spherical harmonics with spherical designs. Finally, we use spherical framelets for denoising Wendland functions as an application, which shows the great potential of spherical designs in spherical data processing.

\end{abstract}

% REQUIRED

% REQUIRED

%\tableofcontents

\section{Introduction}
\label{sec:intro}

Nowadays, the rapid growth of computing power, the massive explosion of data,  and the fast advance of modern information technologies have built a new era of artificial intelligence (AI).  Especially in deep learning, deep neural networks are used to process extremely large and complex data.  Though AI-trained models have many benefits, they require significant computing resources, which can be a challenge for processing certain types of structured data or for certain applications.  For example, data in ImageNet \cite{DengCHI14} are images that can be viewed as samples from a high-dimensional manifold, social network data \cite{goldberg1992using} are with high-dimensional features on a graph-structure domain, data fed into the large language models (LLMs) are typical text type that needs to be vectorized \cite{mikolov2013distributed}.    How to properly represent such data so that one can process them efficiently is the key to their successful applications.  One common way is to use the normalization technique for high-dimensional data.  Once normalized, they can be viewed as data on the $d$-dimensional unit sphere $\mathbb{S}^d:=\{\bm x\in\mathbb{R}^{d+1}\setsep\lVert\bm x\rVert=1\}$, where $\lVert\cdot\rVert$ is the Euclidean norm. Besides these data from normalization,  many data from the real world are actually spherical data, e.g., navigation and geology data, cosmic microwave background data, panoramic images and videos, etc. Therefore, spherical data processing plays an important role in deep learning as well as many other areas \cite{cohen2018spherical,esteves2018learning,kondor2018clebsch,fang2020theory, li2022convolutional}.  

Spherical designs, as one of the very crucial types of point configurations on the sphere, play a very important role in spherical data processing. The increasing demand of processing spherical data has led to increased interest in exploring spherical designs for real-world applications. 
Delsarte, Goethals, and Seidel established the concept of  spherical designs \cite{delsarte1977spherical} in 1977, where a finite point set $X_N:=\{\bm x_1,\ldots,\bm x_N\}\subset\mathbb S^d$ is said to be a {\em spherical $t$-design} if for any polynomial $p:\mathbb{R}^{d+1}\to\mathbb{R}$ of degree at most $t$, the quadrature rule holds:
\begin{equation}
\frac{1}{N}\sum_{i=1}^N p(\bm x_i)=\frac{1}{\omega_d}\int_{\mathbb{S}^d} p(\bm x)\,\mathrm{d}\mu_d(\bm x)\qquad\forall p\in\Pi_t.
\end{equation}
Here, $\mu_d(\bm x)$ denotes the surface measure on $\mathbb{S}^d$ such that $\mu_d(\mathbb S^d):=\omega_d=\frac{2\pi^{\frac{d+1}{2}}}{\Gamma(\frac{d+1}{2})}$ is the surface area of $\mathbb{S}^d$ with $\Gamma(z):=\int_0^\infty x^{z-1}e^{-x}\mathrm{d}x$  and $\Pi_t:=\Pi_t(\mathbb S^d)$ is the space of spherical polynomials on $\mathbb S^d$ with degree at most $t$. As we can see, a spherical $t$-design is an equal weight polynomial-exact quadrature rule in $\Pi_t$ space, which is a  set of points “nicely” distributed on the unit sphere. 
Spherical designs have significant applications in approximation theory, geometry, and combinatorics. For example,  the best packing problems, the minimal energy problems, the optimal configurations related to Smale’s 7th Problem, and so on. 
We refer to a very nice comprehensive survey of Bannai and Bannai \cite{bannai2009survey} for the past five decades of research on spherical designs.  Recently, it has been applied in image reconstruction and signal recovery \cite{chen2018spherical,wang2020tight,xiao2023spherical,zhuang2023denoising}. By leveraging the properties of spherical designs, one can develop efficient algorithms for processing spherical data, which can have a wide range of applications in areas such as computer vision, geophysics, and astrophysics. 

The existence of spherical $t$-designs is a deep theoretical problem that has yielded many profound mathematical results.  Delsarte et al. \cite{delsarte1977spherical} gave the lower bound of a spherical $t$-design on the number of points $N$ for any degree $t\in\mathbb{N}^+$ on $\mathbb{S}^d$: $N\ge \mathcal{O}(t^d)$. 
%\begin{equation*}
%N\ge N^{*}(d,t)=
%\begin{cases}
%2{d+\frac{t-1}{2} \choose d},&\text{$t$ is odd,}\\
%{d+\frac{t}{2} \choose d}+{d+\frac{t}{2}-1 \choose d},&\text{$t$ is even.}
%\end{cases}
%\end{equation*}
%We call a \emph{tight spherical $t$-design} when attains the lower bound. When $d=1$, the vertices of a regular $(t+1)$-gon form a tight spherical $t$-design on the circle $\mathbb S^1$. However, for $d\geq 1$, tight spherical $t$-designs with $N^{*}(d,t)$ points only exist in $t={1,2,3,4,5,7,11}$ \cite{bannai1979tight,bannai1980tight}. 
Seymour and Zaslavsky \cite{seymour1984averaging} proved (non-constructive) that a spherical $t$-design exists for any $t$ if $N$ is sufficiently large. 
Wagner \cite{wagner1991averaging} proposed the feasible upper bounds with $N=\mathcal O(t^{Cd^4})$. Korevaar and Meyers \cite{korevaar1993spherical} provided spherical $t$-designs exist in $N=\mathcal O(t^{d(d+1)/2})$ and conjectured that $N=\mathcal O(t^d)$. Bondarenko, Radchenko and Viazovska \cite{bondarenko2013optimal} verified that spherical $t$-designs indeed exists for $N\geq c_{d}t^{d}$ points, where $c_{d}$ is constant controlled by $d$. Furthermore, they showed that $X_N$ can be {\em well-separated} in the sense that the minimal  separation distance $\delta_{X_N}:=\min_{1\le i<j\le N}\|\bm x_i -\bm x_j\|$  is of order $ \mathcal{O}(N^{-1/d})$ \cite{bondarenko2015well}. 

Besides the theoretical development, numerical algorithms to produce spherical designs are of great importance for spherical data processing.  
%In numerical analysis, Hardin and Sloane \cite{hardin1996mclaren} constructed a sequence of putative spherical $t$-designs with $\frac{1}{2}t^2+o(t^2)$ points. Maier \cite{maier1999numerical} used multiobjective optimization to calculate spherical $t$-designs numerically. Chen and Womersley \cite{chen2006existence} verified the existence of spherical $t$-designs with $(t+1)^2$ points for small $t$. Later, Chen, Frommer and Lang \cite{chen2011computational} showed that spherical $t$-designs with $(t+1)^2$ points exist for all degrees $t$ up to 100. Womersley \cite{womersley2018efficient} constructed symmetric spherical $t$-designs with $N=\frac{t^2+t+4}{2}$ for $t$ up to 325. 
Computing spherical $t$-designs can be regarded as non-linear equations and optimization problems. We refer to \cite{hardin1996mclaren,maier1999numerical,chen2006existence,sloan2009variational,chen2011computational,le2008localized,mhaskar2001spherical,womersley2018efficient,an2020numerical,an2010well,xiao2023spherical} for some of the work on numerical spherical designs. 
In this paper and in what follows, we pay attention to spherical $t$-designs on the most important case $d=2$, i.e., the 2-sphere $\mathbb S^2$. 
Sloan and Womersley \cite{sloan2009variational} presented the variational characterization of spherical $t$-design as a nonnegative quantity $A_{N,t}(X_N)$ given by
\begin{equation}
\label{eq1}
A_{N,t}(X_N):=\frac{4\pi}{N^2}\sum_{\ell=0}^t\sum_{m=-\ell}^{\ell}\left|\sum_{i=1}^N Y_{\ell}^{m}(\bm x_i)\right|^2-1,
\end{equation}
where $Y_\ell^m$ is spherical harmonic with degree $\ell$ and order $m$. They proved important properties about the connection of spherical $t$-designs and $A_{N,t}$. One is that $X_N$ is a spherical $t$-design if and only if $A_{N,t}=0$ (cf. Theorem 3 in \cite{sloan2009variational}), which is the equivalent form with \emph{Weyl sums} \cite{delsarte1977spherical,sloan2009variational} satisfying 
$
\sum_{i=1}^N Y_{\ell}^{m}(\bm x_i)=0,\ \ell= 1,\ldots,t\ \text{and}\ m= -\ell,\ldots,\ell.
$
Hence, finding spherical $t$-designs can be regarded as solving a nonlinear and nonconvex minimization problem:
\begin{equation}
\label{prob}
\min_{X_N\subset\mathbb S^2} A_{N,t}(X_N).
\end{equation}
Based on the quantity $A_{N,t}$, Gr{\"a}f and Potts \cite{graf2011computation} computed numerical spherical $t$-designs using nonequispaced fast spherical Fourier transforms (NFSFTs) and manifold optimization techniques for $t$ up to 1000 with $N \approx\frac{t^2}{2}$ points. Following  \cite{sloan2009variational, graf2011computation}, Xiao and Zhuang \cite{xiao2023spherical} used the trust-region method with NFSFTs to achieve $t$ up to $3200$ numerically.  Numerical results in these papers showed that these methods can approximate spherical $t$-designs with high accuracy.

% main content of this paper.
In this paper, we further investigate the numerical spherical $t$-designs in several aspects.  First, to solve \eqref{prob},  we consider the restart-conjugate gradient method with line search strategy (LS-RCG) in Algorithm~\ref{alg:CGstd} and compare it with the trust-region with preconditioned conjugate gradient method (TR-PCG) presented in \cite{xiao2023spherical}.  Based on algorithms of LS-RCG and TR-PCG , we use two types of point sets on $\mathbb S^2$ as initial point sets for computing spherical $t$-designs and comparing the two algorithms in terms of total iterations, square-root and $\ell_\infty$-error of $A_{N,t}(X_N)$. With the obtained spherical $t$-design point sets, we test their performance in function approximation with functions from the combinations of normalized Wendland functions as well as functions with discontinuity.  Moreover, we further investigate the applications of spherical $t$-designs in the construction of spherical (semi)-tight framelets for spherical image/signal denoising.

% The outline is not required, but we show an example here.
The paper is organized as follows. In Section~\ref{sec:stdmain}, we present the algorithms for computing spherical $t$-designs and provide numerical results for different approaches to compare their performance. In Section~\ref{sec:approx}, we introduce an approximation algorithm using spherical-design point sets and present results for projecting smooth and non-smooth functions onto the sphere. In Section~\ref{sec:apply}, we discuss the application of spherical design point sets in semi-discrete spherical framelet systems and thresholding techniques. We demonstrate how these systems can be used for spherical data processing. The conclusion is given in Section~\ref{sec:conclusions}.

%%%%%%%%%%%%%%%%%%%%%%%%%%%%%%%%%%%%%%
%%%%%%%%%%%%%%%%%%%%%%%%%%%%%%%%%%%%%%
%%%%%%%%%%%%%%%%%%%%%%%%%%%%%%%%%%%%%%
\section{Spherical designs from optimizations}
\label{sec:stdmain}
In this section, we discuss the use of a line-search algorithm with a restart conjugate gradient technique (LS-RCG) for solving the nonlinear and nonconvex minimization problem in \eqref{prob}. We apply such an algorithm to obtain spherical $t$-designs with various initial point sets and compare them with the trust-region method with preconditioned conjugate gradient method (TR-PCG) in \cite{xiao2023spherical}. 

%Before we proceed to the details, we briefly introduce the notation and definition needed for our setting. 

\subsection{Line search with restart conjugate gradient method}
A continuously differentiable function ${f}:\mathbb R^n\to\mathbb R$, the gradient of $f$ at $x$ is defined as
\begin{align}
\label{gd}
\nabla f:=\left[\partial_1f( x),\ldots,\partial_n f( x)\right]^\top,
\end{align}
and the Hessian of $f$ is defined as a $n\times n$ symmetric matrix with elements
\begin{align}
\label{hs}
\left[\mathcal H(f)\right]_{ij}:=\left[\nabla^2 f\right]_{ij}:=\partial_i\partial_jf( x),\qquad 1\leq i,j\leq n,
\end{align}
where $\partial_i$ is the partial derivative with respect to the $i$th coordiante of $x$. 
For the general nonlinear and nonconvex optimization problem:
\begin{equation}
\label{eq:op}
\min_{ x\in \mathbb{R}^d} f( x),
\end{equation}
there are mainly two global convergence approaches: one is the line search, and the other is the trust region.  The trust region with preconditioned conjugate gradient method (TR-PCG) for spherical designs was discussed in \cite{xiao2023spherical}. We refer to \cite{sun2006optimization,hager2006survey,shewchuk1994introduction} for the discussion of trust-region and conjugate gradient methods. In this paper, we use the line search approach with the restart-CG method.  

The iterative scheme using the conjugate gradient (CG) method for solving \eqref{eq:op} is given by
$
 x_{k+1}= x_k+\alpha_k  d_k,
$
where $x_0$ is a starting point, $\alpha_k$ is determined by a line search strategy,  and $d_k$ is a search direction. 
The first seach direction $d_0$ is commonly set as $d_0 = -g_0 :=-\nabla f( x_0)$. Consecutive  $d_k$ is recursively defined as $d_{k+1}=-g_{k+1} +\beta_k d_{k}$ with a scalar $\beta_k$ that can be determined with many choices. The conjugate gradient method was originally proposed by Hestenes and Stiefel \cite{hestenes1952methods} in the 1950s. Without using the restart strategy,  the conjugate gradient method is only linearly convergent \cite{crowder1972linear}.  The restart strategy is used to periodically refresh the algorithm by erasing old information that may be redundant.  It can lead to $n$-step quadratic convergence, that is,
$\|x_{k+n} - x\| = \mathcal{O}(\|x_k - x\|^2)$ \cite{nocedal1999numerical}.
Fletcher and Reeves \cite{fletcher1964function} extended the CG method for unconstrained nonlinear optimization with restart strategy and exact line search, which is the first nonlinear CG method. In 1967, Daniel \cite{daniel1967conjugate} proposed a choice of update parameter $\beta_k$ that required the evaluation of the Hessian in CG method. In this paper, we use $\beta_k = \frac{g_{k+1}^\top B_{k+1} d_k}{d_k^\top B_{k+1} d_k}$ for updating $\beta_k$ with restart criterial $g_k^\top g_{k+1}$ being small, where $B_k$ can be the exact or approximation Hessian, see Aglorithm~\ref{alg:CGstd}. For more discussion on restart procedures for the conjugate gradient method, we refer to \cite{powell1977restart,dai2004restart}.  For $\alpha_k$, we use the commonly used Newton-Raphson method for the line search strategy, see Algorithm~\ref{alg:LSCGstd}.

\renewcommand{\algorithmicrequire}{\textbf{Input:}}
\renewcommand{\algorithmicensure}{\textbf{Output:}}
\begin{algorithm}[htpb!]
  \caption{Line search with restart conjugate gradient method (LS-RCG)}
  \label{alg:CGstd}
  \begin{algorithmic}[1]
    \REQUIRE
      {$x_0$: initial point set;
      $t$: degree;
      $K_{\max}$: maximum iterations;
      $r$: restart orthogonal value;
      $\varepsilon_1$: termination tolerance;

      Initialize $k=0$, $f_0=A_{N,t}(x_0)$, $g_0=\nabla f_0$, $-d_0=g_0$, $B_0\approx\mathcal H(f_0)$.}
    \WHILE{$k\leq K_{\max}$ and $\lVert g_{k}\rVert>\varepsilon_1$}
      \STATE compute step size $\alpha_{k}$ from Algorithm~\ref{alg:LSCGstd}
      \STATE $x_{k+1}=x_k+\alpha_k d_k$, $f_{k+1}=A_{N,t}(x_{k+1})$, $g_{k+1}=\nabla f_{k+1}$ and $B_{k+1}\approx\mathcal H(f_{k+1})$
      \IF{$\lVert g_{k+1}\rVert\leq\varepsilon$}
      \STATE break
      \ENDIF
      \IF {$g_k^\top g_{k+1}\geq r\lVert g_k\rVert^2$}
      \STATE $d_{k+1}=-g_{k+1}$
      \ELSE
      \STATE $\beta_{k}=\frac{g_{k+1}^\top B_{k+1}d_k}{d_k^\top B_{k+1}d_k}$ and $d_{k+1}=-g_{k+1}+\beta_{k}d_{k}$
      \ENDIF
      \STATE $k=k+1$
    \ENDWHILE
    \ENSURE
      {minimizer $x^*\subset\mathbb S^2$.}
  \end{algorithmic}
\end{algorithm}

\begin{algorithm}[htpb!]
  \caption{Line search strategy}
  \label{alg:LSCGstd}
  \begin{algorithmic}[1]
    \REQUIRE
      {$x_k$: current step point set;
      $d_k$: descent direction;
      $I_{\max}$: maximum iterations;
      $\varepsilon_2$: termination tolerance;

      Initialize $n=0$, $d=d_k$, $f_n=A_{N,t}(x_k)$, $g_n=\nabla A_{N,t}(x_k)$, $E_n=\mathcal H(A_{N,t}(x_k))$, $\alpha_n=\frac{g_n^\top d}{d^\top E_n d}$.}
    \WHILE{$n\leq I_{\max}$}
      \STATE $x_{n+1}=x_n+\alpha_n d$ and $f_{n+1}=f(x_{n+1})$
      \IF{$f_{n+1}-f_n>0$}
      \STATE $\alpha_{n+1}=c\alpha_n$, where $c\in (0,1)$
      \ELSE
      \STATE $g_{n+1}=g(x_{n+1})$ and $E_{n+1}=\mathcal H(f_{n+1})$
        \IF{$\frac{\lvert g_{n+1}^\top d \rvert}{\lVert g_{n+1}\rVert \lVert d\rVert}<\varepsilon_2$}
        \STATE break
        \ENDIF
      \STATE $\alpha_{n+1}=\alpha_n-\frac{g_{n+1}^\top d}{d^\top E_{n+1} d}$
      \ENDIF
      \STATE $n=n+1$
    \ENDWHILE
    \ENSURE
      {$\alpha_{n}\in\mathbb R$.}
  \end{algorithmic}
\end{algorithm}

%%%%%%%%%%%%%%%%%%%%%%%%%%%%%%%
%%%%%%%%%%%%%%%%%%%%%%%%%%%%%%% 
\subsection{Spherical harmonics and spherical designs}
For  each spherical coordinate $(\theta,\phi)\in [0,\pi]\times[0,2\pi)$, it is associated with a point  $\bm x:=\bm x(\theta,\phi)=(\sin\theta\cos\phi,\sin\theta\sin\phi,\cos\theta)\in\mathbb{S}^2$. For each $\ell\in\mathbb{N}_0$ and $m=-\ell,\cdots,\ell$, the spherical harmonic $Y_\ell^m: \mathbb{S}^2\rightarrow \C$ can be expressed as
\begin{equation*}
Y_\ell^m(\bm x)=Y_\ell^m(\theta,\phi):=\sqrt{\frac{2\ell+1}{4\pi}\frac{(\ell-m)!}{(\ell+m)!}}P_\ell^m(\cos\theta)\mathrm e^{\mathrm im\phi},
\end{equation*}
where $P_\ell^m:[-1,1]\to\mathbb R$ is the associated Legendre polynomial with form $P_\ell^m(z):=(1-z^2)^{\frac{m}{2}}\frac{\mathrm d^m}{\mathrm d z^m}(P_\ell(z))$ and the Legendre polynomial $P_\ell:[-1,1]\to\mathbb R$ is with form $P_\ell(z)=\frac{1}{2^\ell \ell!}\frac{\mathrm d^\ell}{\mathrm d z^\ell}[(z^2-1)^\ell]$. Note that $Y_0^0=1/\sqrt{4\pi}$. 
The collection $\{Y_\ell^m \setsep \ell\in\mathbb N_0, \lvert m\rvert\leq\ell\}$ forms a complete set of orthonormal basis for the Hilbert space of square-integrable functions $L^2(\mathbb{S}^2):=\{f:\mathbb{S}^2\rightarrow \mathbb{C} \setsep \int_{\mathbb{S}^2} |f(\bm x)|^2 d\mu_2(\bm x)\}$, where the  $L^2$-inner product is defined by $\langle f,g\rangle_{L^2(\mathbb{S}^2)}:=\int_{\mathbb{S}^d} f(\bm x)\overline{g(\bm x)}\mathrm{d}\mu_2(\bm x)$ for $f, g\in L_2(\mathbb{S}^2)$. The orthogonality  $\langle Y_\ell^m,Y_{\ell'}^{m'}\rangle_{L^2(\mathbb{S}^2)}=\delta_{mm'}\delta_{\ell\ell'}$ 
holds for $\ell,\ell' \in \mathbb{N}_0, |m|\le \ell, |m'|\le \ell'$.  Here  $\delta_{ij}$ is Kronecker delta. With such an orthonormal basis, for any $f\in L_2(\mathbb{S}^2)$, it  can be represented as $f=\sum_{\ell=0}^\infty\sum_{m=-\ell}^\ell \hat f_\ell^m Y_\ell^m$ with $\hat f_\ell^m :=\langle f,Y_\ell^m\rangle_{L_2(\mathbb{S}^2)}$.  Moreover,
 $\Pi_t$ is the   linear span of the orthonormal basis $\{Y_\ell^m:\ell\in\mathbb N_0,\ell\leq t,\, \lvert m\rvert\leq\ell\}$ with 
$\text{dim}\,\Pi_t=(t+1)^2$.

In terms of $(\theta, \phi)$, $A_{N,t}(X_N)$ can be regarded as a  function of $2N$ variables. In fact, we can identify the point set $X_N=\{\bm x_1,\ldots, \bm x_N\}\subset \sph^2$ as
\begin{align}
X_N:=(\bm\theta,\bm\phi):=(\theta_1,\ldots,\theta_N,\phi_1,\ldots,\phi_N)
\end{align}
with $\bm\theta = (\theta_1,\ldots,\theta_N)$, $\bm \phi=(\phi_1,\ldots,\phi_N)$, and $\bm x_i:=\bm x_i(\theta_i,\phi_i)$ being the $i$-th point determined by its spherical coordinate satisfying $(\theta_i,\phi_i)\in[0,\pi]\times [0,2\pi)$. The variational characterization $A_{N,t}(X_N)$ in \eqref{eq1} can be written as a smooth function of $2N$ variables:
\begin{align}
\label{ANt:2Nvariables}
A_{N,t}(X_N)&=A_{N,t}(\bm\theta,\bm\phi)%=A_{N,t}(\theta_1,\ldots,\theta_N,\phi_1,\ldots,\phi_N)
%\\&
=\frac{4\pi}{N^2}\sum_{(\ell,m)\in\mathcal{I}_t}\left|\sum_{i\in[N]} Y_{\ell}^{m}(\theta_i,\phi_i)\right|^2-1.
\end{align}
The spherical $t$-design point set $X_N$ can be obtained by solving  \eqref{ANt:2Nvariables} using Algorithm~\ref{alg:CGstd} with the target function $f:=A_{N,t}(\bm\theta,\bm\phi)$. As discussed in \cite{xiao2023spherical}, the key is the fast evaluations of $A_{N,t}$, $\nabla A_{N,t}$, and $\mathcal{H}(A_{N,t})$. By \cite[Theorems 2.2 and 2.3]{xiao2023spherical}, they can be evaluated through the fast spherical harmonic transforms with arithmetic complexity of order $\mathcal{O}(t^2\log^2{t}+N\log^2(\frac{1}{\epsilon}))$, where $\epsilon$ is prescribed accuracy of the approximate algorithms \cite{kunis2003fast,potts2001fast,ware1998fast,keiner2008fast}. 

Moreover, the Hessian $\mathcal{H}(A_{N,t}) =\mathcal{H}_1+\mathcal{H}_2$   is the  sum of a diagonal matrix $\mathcal{H}_1$ and a rank-one matrix $\mathcal{H}_2$, see \cite[Theorem 2.3]{xiao2023spherical}.  In Algorithm~\ref{alg:CGstd}, $B_k$ is not necessary to be the exact full Hessian. One can use an approximation version of the Hessian for $B_k$, e.g., $\mathcal{H}_2$.

%%%%%%%%%%%%%%%%%%%%%%%%%%%%%%%%%%%%%%%%%%%%%
\subsection{Numerical spherical designs from two different approaches}
In this section, based on LS-RCG in Algorithm~\ref{alg:CGstd} and TR-PCG from \cite[Algorithm 2.1]{xiao2023spherical} to solve \eqref{prob}, we get numerical spherical $t$-designs\footnote{All numerical experiments in this paper are conducted in MATLAB R2023a on a macOS Sonoma system's iMac (Retina 5K, 27-inch, 2020) with Intel Core i5 10500 CPU and 32 GB DDR4 memory.}. We use two types of point sets on $\mathbb S^2$ as initial point sets for computing spherical $t$-designs with the setting as follows.
\begin{enumerate}[I.]
\item Spiral points (SP). To generate spiral points $\bm x_n=(\theta_n,\phi_n)$ on $\mathbb S^2$ for $n=1,\ldots,N$, we set
\begin{equation*}
\theta_n :=\arccos\left(\frac{2n-(N+1)}{N}\right),\\
\phi_n :=\pi(2n-(N+1))\varphi^{-1},
\end{equation*}
where $\varphi=\frac{1+\sqrt{5}}{2}$ is the golden ratio, refer to \cite{swinbank2006fibonacci} which is Fibonacci spiral on the sphere, same on the setting of the initial distribution of spiral points in \cite{graf2011computation}.
\item Uniformly distributed points (UD). To generate uniformly distributed points on unit sphere $\mathbb S^2$, we need to make sure that for each surface area $\mathrm{d}\mu=\sin\theta\mathrm{d}\theta\mathrm{d}\phi$ contain the same number of points. Thus by \cite{weisstein2004Uniformrand}, we generate random variables $k_n\in (0,1)$ and $p_n\in (0,1)$ for $n=1,\ldots,N$, then we have
    \begin{equation*}
    \theta_n :=\arccos\left(1-2k_n \right),\\
    \phi_n :=2\pi p_n.
    \end{equation*}
\end{enumerate}
We deal with point set $X_N\subset\mathbb S^2$ from the beginning by fixing the first point $\bm x_1=(0,0)\in X_N$ be the north pole point and the second point $\bm x_2=(\theta_2,0)\in X_N$ be on the prime meridian. Then  we let $X_N =(\theta_2,\ldots,\theta_N,\phi_3,\ldots,\phi_N)^\top=[\Theta_{N-1},  \Phi_{N-2}]^\top$. 

We show in Tables~\ref{table1} and \ref{table2} for various degree $t$  with $N=(t+1)^2$,  including the performance of different methods with the full  Hessian $\mathcal{H}(A_{N,t})=\mathcal{H}_1+\mathcal{H}_2$ and the approximation Hessian $\mathcal{H}_2$ (the rank-one matrix) under different point sets. We compare the results of different initial point sets and the efficiency and accuracy of different approaches in the two methods.  

In terms of the accuracy of function values $A_{N,t}(X_N)$ and $\nabla A_{N,t}(X_N)$ in the two algorithms, TR-PCG and LS-RCG have similar performance (both up to order 12 accuracy). In terms of the total iteration $K$, TR-PCG uses significantly fewer iterations than those of LS-RCG.  The LS-RCG uses a lot of iterations in the line-search step while TR-PCG's main time-consuming step is on the step of finding the trust regions. In terms of the accuracy of using full Hessian and approximation Hessian, there are not many differences between these two approaches.  However, in both methods, approximation Hessian takes fewer steps compared to full Hessian. Generally, full Hessian is a bit better than the approximation Hessian in accuracy.  In terms of the point sets, from both tables and Figures, one can see that structured initial point sets SPs use significantly less iterations than those of the randomly distributed initial point sets UDs.

\begin{table}[htpb!]
\centering
\caption{Computing of spherical $t$-designs by different approaches from SPD point sets. In point sets $X_N$, SPD means generated by SP. $K$ is the total iterations in Algorithms (Alg). $\mathcal H$ represents the type of Hessian (full or approximation).}\label{table1}
\begin{small}
\begin{tabular}{lll|lllcc}
\hline														
$X_N$  &  Alg  &  $\mathcal H$  &  $t$  &  $N$  &  $K$  &  $\sqrt{A_{N,t}(X_N)}$  &  $\lVert\nabla A_{N,t}(X_N)\rVert_{\infty}$  \\
\hline              
\multirow{10}{*}{SPD}  &  \multirow{10}{*}{LS-RCG}  & \multirow{10}{*}{Full}  & 10 & 121 & 3551 & 3.95E-12 &  1.79E-14 \\ 
~  &  ~  &  ~  & 20 & 441 & 4111 & 1.64E-12 &  1.25E-14 \\ 
 ~  &  ~  &  ~  & 30 & 961 & 5484 & 1.49E-12 &  7.01E-15 \\ 
~  &  ~  &  ~  & 40 & 1681 & 7360 & 1.35E-12 &  6.77E-15 \\ 
~  &  ~  &  ~  & 50 & 2601 & 10397 & 1.60E-12 &  6.34E-15 \\ 
~  &  ~  &  ~  & 60 & 3721 & 11663 & 1.09E-12 &  4.27E-15 \\ 
~  &  ~  &  ~  & 70 & 5041 & 13409 & 1.39E-12 &  4.50E-15 \\ 
~  &  ~  &  ~  & 80 & 6561 & 13709 & 1.55E-12 &  5.63E-15 \\ 
~  &  ~  &  ~  & 90 & 8281 & 16380 & 1.27E-12 &  4.11E-15 \\ 
~  &  ~  &  ~  & 100 & 10201 & 18866 & 1.01E-12 &  3.11E-15 \\ 
\hline              
\multirow{10}{*}{SPD}  &  \multirow{10}{*}{LS-RCG}  &  \multirow{10}{*}{Appr}  & 10 & 121 & 3652 & 3.95E-12 &  1.87E-14 \\ 
~  &  ~  &  ~  & 20 & 441 & 4180 & 1.64E-12 &  1.19E-14 \\ 
~  &  ~  &  ~  & 30 & 961 & 5336 & 1.49E-12 &  7.08E-15 \\ 
~  &  ~  &  ~  & 40 & 1681 & 7270 & 1.35E-12 &  4.99E-15 \\ 
~  &  ~  &  ~  & 50 & 2601 & 10399 & 1.60E-12 &  6.49E-15 \\ 
~  &  ~  &  ~  & 60 & 3721 & 11670 & 1.09E-12 &  4.38E-15 \\ 
~  &  ~  &  ~  & 70 & 5041 & 13486 & 1.39E-12 &  4.15E-15 \\ 
~  &  ~  &  ~  & 80 & 6561 & 13419 & 1.55E-12 &  2.64E-15 \\ 
~  &  ~  &  ~  & 90 & 8281 & 16546 & 1.27E-12 &  2.59E-15 \\ 
~  &  ~  &  ~  & 100 & 10201 & 18867 & 1.01E-12 &  3.53E-15 \\ 
\hline              
\hline
\multirow{10}{*}{SPD}  &  \multirow{10}{*}{TR-PCG}  &  \multirow{10}{*}{Full}  & 10 & 121 & 169 & 3.94E-12 &  2.44E-14 \\ 
~  &  ~  &  ~  & 20 & 441 & 379 & 1.71E-12 &  9.44E-16 \\ 
~  &  ~  &  ~  & 30 & 961 & 480 & 1.49E-12 &  9.64E-16 \\ 
~  &  ~  &  ~  & 40 & 1681 & 727 & 1.34E-12 &  8.77E-16 \\ 
~  &  ~  &  ~  & 50 & 2601 & 764 & 1.58E-12 &  9.39E-15 \\ 
~  &  ~  &  ~  & 60 & 3721 & 939 & 1.08E-12 &  2.67E-15 \\ 
~  &  ~  &  ~  & 70 & 5041 & 1355 & 1.39E-12 &  9.78E-16 \\ 
~  &  ~  &  ~  & 80 & 6561 & 1235 & 1.54E-12 &  1.13E-15 \\ 
~  &  ~  &  ~  & 90 & 8281 & 1371 & 1.26E-12 &  6.89E-16 \\ 
~  &  ~  &  ~  & 100 & 10201 & 1699 & 1.00E-12 &  8.51E-16 \\ 
\hline              
\multirow{10}{*}{SPD}  &  \multirow{10}{*}{TR-PCG}  &  \multirow{10}{*}{Appr}  & 10 & 121 & 200 & 3.95E-12 &  7.77E-16 \\ 
~  &  ~  &  ~  & 20 & 441 & 377 & 1.65E-12 &  6.03E-16 \\ 
~  &  ~  &  ~  & 30 & 961 & 458 & 1.49E-12 &  9.76E-16 \\ 
~  &  ~  &  ~  & 40 & 1681 & 649 & 1.35E-12 &  9.74E-16 \\ 
~  &  ~  &  ~  & 50 & 2601 & 743 & 1.58E-12 &  1.07E-14 \\ 
~  &  ~  &  ~  & 60 & 3721 & 906 & 1.08E-12 &  2.83E-15 \\ 
~  &  ~  &  ~  & 70 & 5041 & 1037 & 1.39E-12 &  1.37E-15 \\ 
~  &  ~  &  ~  & 80 & 6561 & 1121 & 1.54E-12 &  6.45E-16 \\ 
~  &  ~  &  ~  & 90 & 8281 & 1385 & 1.26E-12 &  8.76E-16 \\ 
~  &  ~  &  ~  & 100 & 10201 & 1461 & 1.00E-12 &  9.97E-16 \\ 
\hline
\end{tabular}
\end{small}
\end{table}

 \begin{table}[htpb!]
\centering
\caption{Computing of spherical $t$-designs by different approaches from SUD point set. In point sets $X_N$, SUD means generated by UD. $K$ is the total iterations in Algorithms (Alg). $\mathcal H$ represents the type of Hessian (full or approximation).}\label{table2}
\begin{small}
\begin{tabular}{lll|lllcc}
\hline														
$X_N$  &  Alg  &  $\mathcal H$  &  $t$  &  $N$  &  $K$  &  $\sqrt{A_{N,t}(X_N)}$  &  $\lVert\nabla A_{N,t}(X_N)\rVert_{\infty}$  \\
\hline                         
\multirow{5}{*}{SUD}  &  \multirow{5}{*}{LS-RCG}  &  \multirow{5}{*}{Full}  & 10 & 121 & 3633 & 3.38E-12 &  1.85E-14 \\ 
~  &  ~  &  ~  & 20 & 441 & 5167 & 1.96E-12 &  9.65E-15 \\ 
~  &  ~  &  ~  & 30 & 961 & 7232 & 1.79E-12 &  6.88E-15 \\ 
~  &  ~  &  ~  & 40 & 1681 & 8814 & 1.54E-12 &  9.50E-15 \\ 
~  &  ~  &  ~  & 50 & 2601 & 20221 & 1.64E-12 &  5.85E-15 \\ 
\hline              
\multirow{5}{*}{SUD}  &  \multirow{5}{*}{LS-RCG}  &  \multirow{5}{*}{Appr}  & 10 & 121 & 3701 & 4.28E-12 &  1.96E-14 \\ 
~  &  ~  &  ~  & 20 & 441 & 5233 & 1.79E-12 &  1.32E-14 \\ 
~  &  ~  &  ~  & 30 & 961 & 6595 & 1.81E-12 &  6.80E-15 \\ 
~  &  ~  &  ~  & 40 & 1681 & 8811 & 1.60E-12 &  4.86E-15 \\ 
~  &  ~  &  ~  & 50 & 2601 & 14260 & 1.57E-12 &  5.92E-15 \\ 
\hline              
\hline
\multirow{5}{*}{SUD}  &  \multirow{5}{*}{TR-PCG}  &  \multirow{5}{*}{Full}  & 10 & 121 & 316 & 4.16E-12 &  2.92E-15 \\ 
~  &  ~  &  ~  & 20 & 441 & 794 & 1.98E-12 &  9.15E-16 \\ 
~  &  ~  &  ~  & 30 & 961 & 847 & 1.72E-12 &  3.17E-15 \\ 
~  &  ~  &  ~  & 40 & 1681 & 1546 & 1.57E-12 &  2.57E-15 \\ 
~  &  ~  &  ~  & 50 & 2601 & 1660 & 1.44E-12 &  1.74E-14 \\ 
\hline              
\multirow{5}{*}{SUD}  &  \multirow{5}{*}{TR-PCG}  &  \multirow{5}{*}{Appr}  & 10 & 121 & 234 & 3.87E-12 &  7.54E-16 \\ 
~  &  ~  &  ~  & 20 & 441 & 567 & 1.80E-12 &  1.19E-15 \\ 
~  &  ~  &  ~  & 30 & 961 & 765 & 1.59E-12 &  3.18E-15 \\ 
~  &  ~  &  ~  & 40 & 1681 & 849 & 1.61E-12 &  1.91E-15 \\ 
~  &  ~  &  ~  & 50 & 2601 & 1151 & 1.47E-12 &  2.71E-14 \\ 
\hline														
\end{tabular}
\end{small}
\end{table}

\begin{figure}[htpb!]
\centering
\subfigure[SP point set 3D. Left: Initial. Right: Final]{
\label{Fig.sub.sp}
\includegraphics[width=0.44\textwidth]{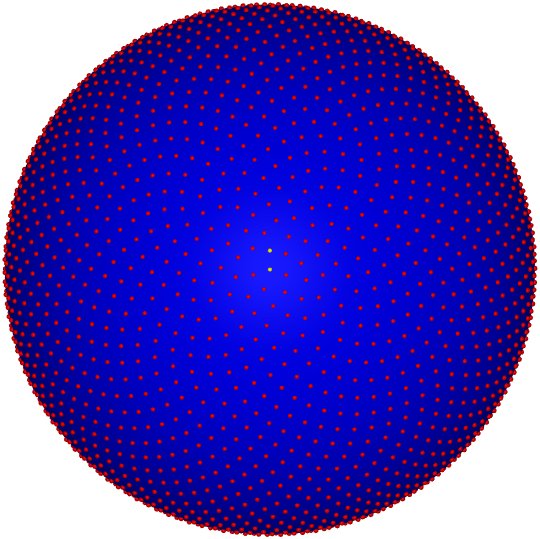}\quad
\includegraphics[width=0.44\textwidth]{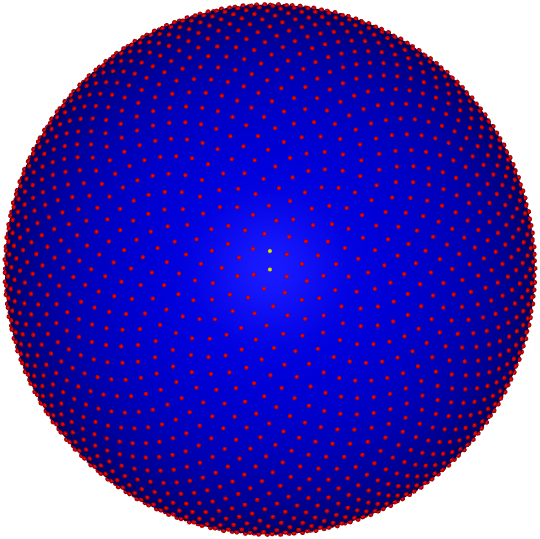}
}
\\
\subfigure[SP point set 2D. Left: Initial. Right: Final]{
\label{Fig.sub.sp2}
\includegraphics[width=0.44\textwidth]{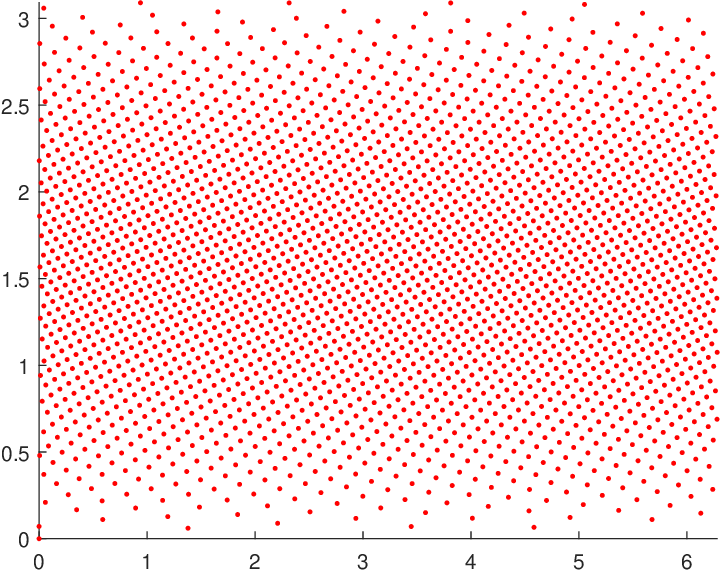}\quad
\includegraphics[width=0.44\textwidth]{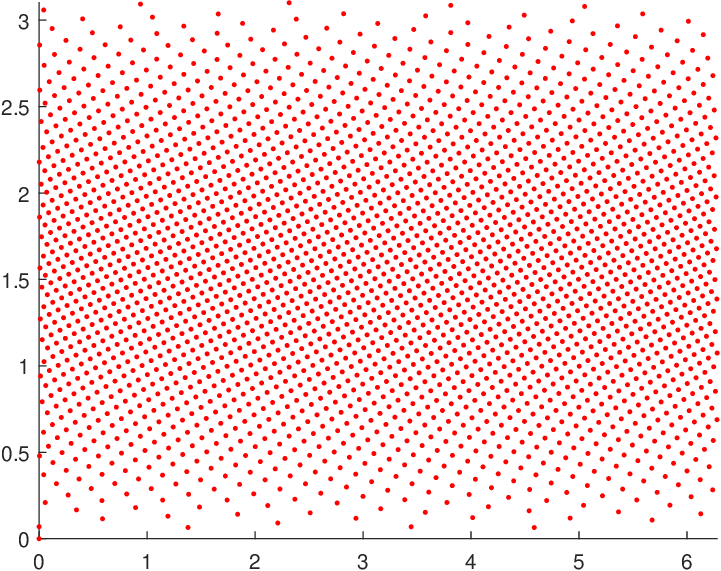}
}
%\\
%\subfigure[UD point set 3D. Left: Initial. Right: Final]{
%\label{Fig.sub.ud}
%\includegraphics[width=0.4\textwidth]{50_2601_UD_initial.png}
%\includegraphics[width=0.4\textwidth]{50_2601_UD_final.png}
%}
%\\
%\subfigure[UD point set 2D. Left: Initial. Right: Final]{
%\label{Fig.sub.ud2}
%\includegraphics[width=0.44\textwidth]{50_2601_UD_initial_2D.eps}
%\includegraphics[width=0.44\textwidth]{50_2601_UD_final_2D.eps}
%}

\caption{Numerical simulation of spherical initial point sets (left column) as input point sets and related spherical $t$-design point sets (right column) as output point sets based on TR-PCG for $t=50$ and $N=(t+1)^2$ on $\mathbb S^2$. (a) on sphere. (b) equirectangular projection. }
\label{Fig.label.ptsp}
\end{figure}

\begin{figure}[htpb!]
\centering
%\subfigure[SP point set 3D. Left: Initial. Right: Final]{
%\label{Fig.sub.sp}
%\includegraphics[width=0.4\textwidth]{50_2601_SP_initial.png}
%\includegraphics[width=0.4\textwidth]{50_2601_SP_final.png}
%}
%\\
%\subfigure[SP point set 2D. Left: Initial. Right: Final]{
%\label{Fig.sub.sp2}
%\includegraphics[width=0.44\textwidth]{50_2601_SP_initial_2D.eps}
%\includegraphics[width=0.44\textwidth]{50_2601_SP_final_2D.eps}
%}
%\\
\subfigure[UD point set 3D. Left: Initial. Right: Final]{
\label{Fig.sub.ud}
\includegraphics[width=0.44\textwidth]{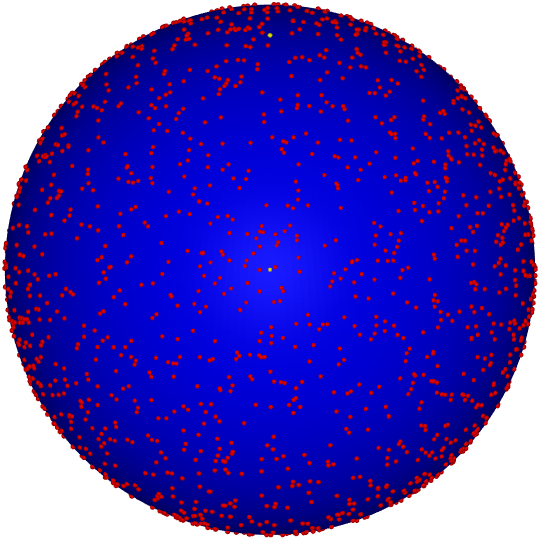}\quad
\includegraphics[width=0.44\textwidth]{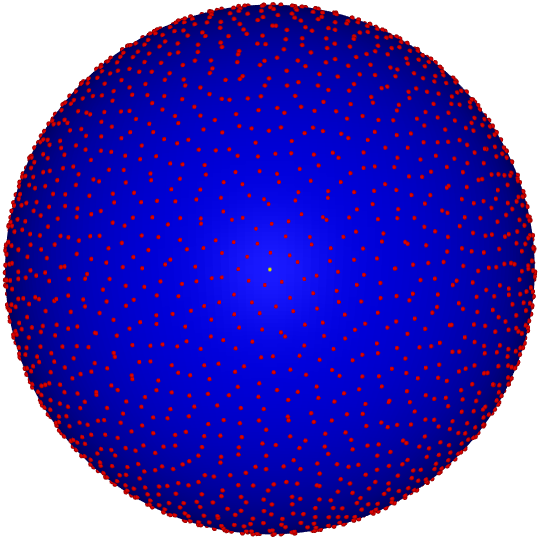}
}
\\
\subfigure[UD point set 2D. Left: Initial. Right: Final]{
\label{Fig.sub.ud2}
\includegraphics[width=0.44\textwidth]{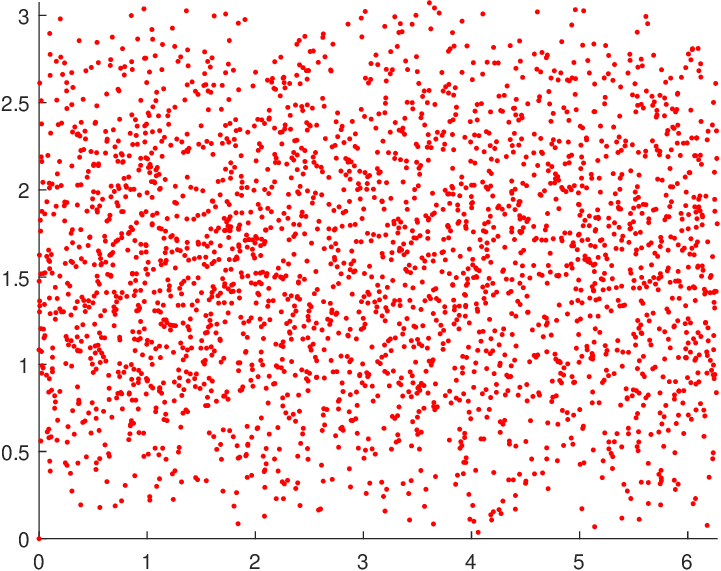}\quad
\includegraphics[width=0.44\textwidth]{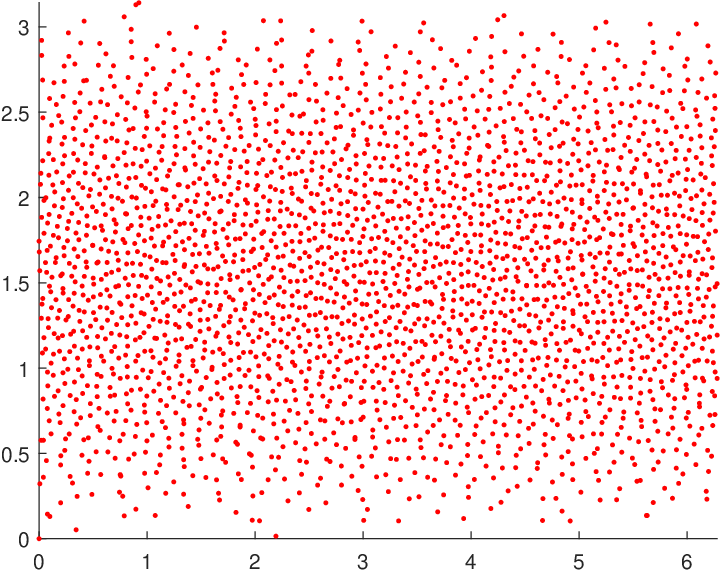}
}

\caption{Numerical simulation of spherical initial point sets (left column) as input point sets and resulted spherical $t$-design point sets (right column) as output point sets based on TR-PCG for $t=50$ and $N=(t+1)^2$ on $\mathbb S^2$. (a): on sphere. (b): equirectangular projection. }
\label{Fig.label.ptud}
\end{figure}

Based on these tables, we further show the initial point sets and related spherical $t$-design point sets on Figures~\ref{Fig.label.ptsp} and ~\ref{Fig.label.ptud} by choosing $t=50$ and $N=(t+1)^2$ for SP and UD. One can see that the SP initial point set is well distributed already and the pattern of the final spherical design does not change much from such an SP initial point set. On the other hand, one can see that the final spherical design of the UD point set is well distributed compared to its initial stage.

%%%%%%%%%%%%%%%%%%%%%%%%%%%%%%%%%%%%%%
%%%%%%%%%%%%%%%%%%%%%%%%%%%%%%%%%%%%%%
%%%%%%%%%%%%%%%%%%%%%%%%%%%%%%%%%%%%%%
\section{Function approximation with spherical designs}
\label{sec:approx}
In this section, we study the approximation of function using polynomials with the obtained spherical designs. 

\subsection{Polynomial spaces projection}
\label{subsec:proj}
First, we consider the orthogonal projection on $\Pi_t$. Given a  spherical data $\{(x_i,f(x_i)) \setsep  x_i\in X_N\}$ from some spherical signal  $f: \mathbb{S}^2 \rightarrow \mathbb{R}$ in $L_2(\mathbb{S}^2)$ and a spherical $t$-design $X_N$ obtained above,  the spherical signal $f$ can be projected on $\Pi_t$   through orthogonal projection. That is, $f=f_0+f_1$ as 
\begin{align}
f(x_i)=\sum_{\ell=0}^t\sum_{m=-\ell}^\ell \hat c_\ell^m Y_\ell^m(\theta_i,\phi_i)+f_1(x_i),
\end{align}
where $f_0(x)=\sum\limits_{\ell=0}^t\sum\limits_{m=-\ell}^\ell \hat c_\ell^m Y_\ell^m(\theta,\phi)\in\Pi_t$, and the residual $f_1\notin\Pi_t$.  Here $\hat f_\ell^m:=\langle f, Y_\ell^m\rangle$. To obtain $f_0$, i.e., the coefficient set $\hat {\mathbf c}=\{\hat c_\ell^m: \ell=0,\ldots, t, |m|\le \ell\}$, we can use the least squares method through the minimization problem
\begin{equation}
\label{ls}
\min_{\mathbf f_0\in\Pi_t}\lVert \mathbf f-\mathbf f_0\rVert,
\end{equation}
where $\mathbf{f}:=f|_{X_N}$ and $\mathbf{f}_0:=f_0|_{X_N}$ are the vector representations of $f$ and $f_0$ on $X_N$, respectively. 
To solve \eqref{ls}, by using matrix representation $\mathbf f_0=\mathbf Y_t\hat{\mathbf c}$, where $\mathbf Y_t$ is the $N\times (t+1)^2$ matrix determined by $(Y_\ell^m(x_i))_{x_i\in X_N; |\ell|\le t, |m|\le \ell}$,  the least squares problem  in \eqref{ls} can be solved through the following equation:
\begin{align}
\label{eq20}
\mathbf Y_t^*(\mathbf w\odot\mathbf Y_t\hat{\mathbf c})=\mathbf Y_t^*(\mathbf w\odot\mathbf f)
\end{align}
with a wight $\mathbf w$. It can be done by using the CG method for the equation $Ax=b$ with 
 $A=\mathbf Y_t^* \mathrm{diag}(\mathbf w)\mathbf Y_t$,  $x=\hat{\mathbf c}$ and $b=\mathbf Y_t^*(\mathbf{w}\odot \mathbf f$). See Algorithm~\ref{alg:projCG}. For convenience, we acquiesce weight function $W(X_N):=\mathbf w$ in corresponding operation $(\mathbf Y_t^*\cdot)$ on arbitrary point set $X_N$ in Algorithm~\ref{alg:projCG}.

\begin{algorithm}[htpb!]
  \caption{Projection by Conjugate Gradient Algorithm}
  \label{alg:projCG}
  \begin{algorithmic}[1]
    \REQUIRE
      {$t$: polynomial degree;
      $X_N$: spherical point set;
      $W$: weight function of $X_N$;
      $K_{\max}$: maximum iterations;
      $\varepsilon$: termination tolerance;

      Initialize $x=0$, $k=0$, $r_0=\mathbf Y_t^*(\mathbf w\odot \mathbf f)$, $A= \mathbf Y_t \mathrm{diag}(\mathbf w)\mathbf Y_t$.}
    \WHILE{$\lVert r_{k+1}\rVert>\varepsilon$ and $k\leq K_{\max}$}
      \IF{$k=0$}
      \STATE $p_1=r_0$
      \ELSE
      \STATE $p_{k+1}=r_k+\frac{\lVert r_k\rVert^2}{\lVert r_{k-1}\rVert^2}p_k$
      \ENDIF
      \STATE compute step size $\alpha=\frac{\lVert r_k\rVert^2}{p_{k+1}^\top Ap_{k+1}}$
      \STATE $x_{k+1}=x_k+\alpha p_{k+1}$ and $r_{k+1}=r_k-\alpha Ap_{k+1}$
      \STATE $k=k+1$
    \ENDWHILE
    \ENSURE
      {$x^*\in\mathbb R^n$.}
  \end{algorithmic}
\end{algorithm}

%%%%%%%%%%%%%%%%%%%%%%%%%%%%%%%%%%%%%%
%%%%%%%%%%%%%%%%%%%%%%%%%%%%%%%%%%%%%%
\subsection{Approximation of spherical functions from Wendland functions}
\label{subsec:smth}
Next, we consider the approximation of some spherical functions using polynomials with different point sets. The functions are the combinations of normalized Wendland functions from a family of compactly supported radial basis functions (RBF). Let $(t)_{+}:=\max\{t,0\}$ for $t\in\mathbb{R}$. The original Wendland functions are
\begin{equation*}
  \tilde\phi_k(t) := \begin{cases}
  (1-t)_{+}^{2}, & k = 0,\\[1mm]
  (1-t)_{+}^{4}(4t + 1), & k = 1,\\[1mm]
  \displaystyle (1-t)_{+}^{6}(35t^2 + 18t + 3)/3, & k = 2,\\[1mm]
  (1-t)_{+}^{8}(32t^3 + 25t^2 + 8t + 1), & k = 3,\\[1mm]
  \displaystyle (1-t)_{+}^{10}(429t^4 + 450t^3 + 210t^2 + 50t + 5)/5, & k = 4.
  \end{cases}
\end{equation*}
The normalized (equal area) Wendland functions are
\begin{equation*}
    \phi_k(t) := \tilde\phi_k\Bigl(\frac{t}{\delta_{k}}\Bigr),\quad \delta_{k} := \frac{(3k+3)\Gamma(k+\frac{1}{2})}{2\:\Gamma(k+1)},\quad k\geq0.
\end{equation*}
The Wendland functions $\phi_k(t)$ pointwise converge to Guassian when $k\to\infty$, refer to \cite{chernih2014wendland}. Thus the main change as $k$ increases is the smoothness of $f$. Let $\bm z_1:=(1,0,0),\bm z_2:=(-1,0,0),\bm z_3:=(0,1,0),\bm z_4:=(0,-1,0),\bm z_5:=(0,0,1),\bm z_6:=(0,0,-1)$ be regular octahedron vertices and define \cite{gia2010multiscale}
\begin{equation}\label{eq:Phi}
 f_{k}(\bm x)
  := \sum_{i=1}^{6}\phi_k(\|\bm z_{i} - \bm x\|), \quad k\geq0.
\end{equation}
 The paper \cite{gia2010multiscale} has proved that $f_{k}\in H^{k+\frac{3}{2}}(\mathbb S^2)$, where $H^{\sigma}(\mathbb S^2):=\{f\in L_2(\mathbb S^2):\sum\limits_{\ell=0}^\infty\sum\limits_{\lvert m\rvert\leq\ell} (1+\ell)^{2\sigma}\lvert\hat f_{\ell}^m\rvert^2\}<\infty\}$ is the Sobolev space with smooth parameter $\sigma>1$. The function $f_k$ has limited smoothness at the centers $\bm z_i$ and at the boundary of each cap with center $\bm z_i$. These features make $f_{k}$ relatively difficult to approximate in these regions, especially for small $k$.

Given a point set $X_N$ with weight $\mathbf w$ (not necessary quadrature point sets), we can sample $f_k$ on  $X_N$ as spherical signal $v: X_N\rightarrow \mathbb{R}$, and compute the projection $v_1$ and residual $r=v-v_1$ based on Algorithm~\ref{alg:projCG}, where we set the maximum iterations $K_{\max}=1000$ and termination tolerance $\varepsilon$ = 2.2204e-16 (floating-point relative accuracy of MATLAB). To check whether the point sets have good approximation properties for such functions, we use the relative projection $L_2$-error defined by
\begin{equation*}
err(v_1,v)=\frac{\lVert v-v_1\rVert}{\lVert v\rVert}=\frac{\lVert r\rVert}{\lVert v\rVert}.
\end{equation*}

We show the results in Table~\ref{table:err2} under the setting of input polynomial degree $T=\frac{t}{2}$ and weight $W=\mathbf w$ on Algorithm~\ref{alg:projCG}, and set input polynomial degree $T=\frac{t}{2}$ and weight $W=\sqrt{\mathbf w}$ on Algorithm~\ref{alg:projCG} to obtain Table~\ref{table:err4}. For comparing the initial points, we set SP and UD with equal weight $\mathbf w \equiv \frac{4\pi}{N}$. The spherical design point sets are equal-weight quadrature rules. Also in Tables~\ref{table:err2} and \ref{table:err4}, we give the results of spherical $t$-designs with degree $t\approx 200,400$ and $N\approx (t+1)^2$ in different initial points: Spiral points (SP), uniform distributed points (UD), Icosahedron vertices mesh points (IV), and HEALpix points (HL). ((SP, SPD), (UD, SUD), (IV, SID), and (HL, SHD)  are the initial point set and final spherical $t$-design point set via TR-PCG, respectively).  SID means generated by IV, and SHD means generated by HL) For more details on these point sets, we refer to \cite{xiao2023spherical}. We shall mention that the paper \cite{xiao2023spherical} has shown the partial results of $t\approx 200$ in Table~\ref{table:err2}.

\begin{table}[htpb!]
    \centering
    \caption{Relative $L_{2}$-errors $err(f_T,f_k)$ for Wendland functions $f_0,\ldots,f_4$ approximated by $\Pi_T$ functions under the setting $T=\frac{t}{2}$ and $W=\mathbf w$ on Algorithm~\ref{alg:projCG}.}\label{table:err2}
    \begin{small}
    \begin{tabular}{llllllll}
    \hline
        $t$ & $N$ & $Q_N$ & $f_0$ & $f_1$ & $f_2$ & $f_3$ & $f_4$ \\
        \hline
        200 & 40401 & SP & 5.64E-04 & 3.19E-06 & 5.25E-08 & 3.39E-09 & 3.21E-09 \\
        200 & 40401 & SPD & 5.78E-04 & 3.20E-06 & 5.25E-08 & 1.69E-09 & 8.92E-11 \\
        400 & 160801 & SP & 1.49E-04 & 2.05E-07 & 1.47E-09 & 1.33E-09 & 1.45E-09 \\ 
        400 & 160801 & SPD & 1.49E-04 & 2.06E-07 & 8.61E-10 & 7.39E-12 & 2.28E-12 \\ 
        \hline
        200 & 40401 & UD & 6.09E-04 & 3.07E-06 & 8.50E-08 & 8.71E-08 & 8.53E-08 \\
        200 & 40401 & SUD & 6.99E-04 & 3.51E-06 & 5.63E-08 & 1.79E-09 & 1.07E-10 \\
        400 & 160801 & UD & 1.22E-04 & 1.86E-07 & 8.07E-08 & 8.55E-08 & 9.02E-08 \\ 
        400 & 160801 & SUD & 1.68E-04 & 2.19E-07 & 1.23E-09 & 9.33E-10 & 1.02E-09 \\ 
        \hline
        201 & 40962 & IV & 8.12E-04 & 3.32E-06 & 5.28E-08 & 4.57E-09 & 6.10E-08 \\
        201 & 40962 & SID & 6.15E-04 & 3.11E-06 & 5.08E-08 & 1.64E-09 & 8.74E-11 \\
        403 & 163842 & IV & 2.05E-04 & 2.14E-07 & 3.98E-09 & 3.88E-09 & 5.59E-08 \\ 
        403 & 163842 & SID & 1.80E-04 & 2.18E-07 & 8.81E-10 & 7.46E-12 & 2.29E-12 \\ 
        \hline
        220 & 49152 & HL & 5.98E-04 & 2.28E-06 & 3.18E-08 & 8.68E-09 & 8.28E-09 \\
        220 & 49152 & SHD & 5.98E-04 & 2.28E-06 & 3.04E-08 & 8.11E-10 & 3.59E-11 \\
        442 & 196608 & HL & 1.46E-04 & 1.45E-07 & 8.19E-10 & 3.94E-10 & 3.11E-10 \\ 
        442 & 196608 & SHD & 1.46E-04 & 1.45E-07 & 4.89E-10 & 4.07E-12 & 2.39E-12 \\ 
        \hline
    \end{tabular}
    \end{small}
\end{table}

\begin{table}[htpb!]
    \centering
    \caption{Relative $L_{2}$-errors $err(f_T,f_k)$ for Wendland functions $f_0,\ldots,f_4$ approximated by $\Pi_T$ functions under the setting $T=\frac{t}{2}$ and $W=\sqrt{\mathbf w}$ on Algorithm~\ref{alg:projCG}.}\label{table:err4}
    \begin{small}
    \begin{tabular}{llllllll}
    \hline
        $t$ & $N$ & $Q_N$ & $f_0$ & $f_1$ & $f_2$ & $f_3$ & $f_4$ \\
        \hline
        200 & 40401 & SP & 5.64E-04 & 3.19E-06 & 5.24E-08 & 1.69E-09 & 9.34E-11 \\
        200 & 40401 & SPD & 5.78E-04 & 3.20E-06 & 5.25E-08 & 1.69E-09 & 8.92E-11 \\
        400 & 160801 & SP & 1.49E-04 & 2.05E-07 & 8.58E-10 & 1.41E-11 & 1.32E-11 \\ 
        400 & 160801 & SPD & 1.49E-04 & 2.06E-07 & 8.61E-10 & 7.39E-12 & 2.25E-12 \\ 
        \hline
        200 & 40401 & UD & 6.09E-04 & 3.07E-06 & 4.80E-08 & 1.89E-09 & 1.46E-09 \\
        200 & 40401 & SUD & 6.99E-04 & 3.51E-06 & 5.63E-08 & 1.79E-09 & 1.07E-10 \\
        400 & 160801 & UD & 1.22E-04 & 1.76E-07 & 9.04E-10 & 6.74E-10 & 6.49E-10 \\ 
        400 & 160801 & SUD & 1.68E-04 & 2.19E-07 & 8.95E-10 & 7.35E-12 & 1.34E-12 \\ 
        \hline
        201 & 40962 & IV & 8.12E-04 & 3.32E-06 & 5.26E-08 & 1.71E-09 & 2.95E-10 \\
        201 & 40962 & SID & 6.15E-04 & 3.11E-06 & 5.08E-08 & 1.64E-09 & 8.74E-11 \\
        403 & 163842 & IV & 2.05E-04 & 2.14E-07 & 9.01E-10 & 2.45E-10 & 2.47E-10 \\ 
        403 & 163842 & SID & 1.80E-04 & 2.18E-07 & 8.81E-10 & 7.46E-12 & 2.29E-12 \\ 
        \hline
        220 & 49152 & HL & 5.98E-04 & 2.28E-06 & 3.04E-08 & 8.12E-10 & 5.82E-11 \\
        220 & 49152 & SHD & 5.97E-04 & 2.28E-06 & 3.04E-08 & 8.11E-10 & 3.59E-11 \\
        442 & 196608 & HL & 1.46E-04 & 1.45E-07 & 4.89E-10 & 3.94E-10 & 3.11E-10 \\ 
        442 & 196608 & SHD & 1.46E-04 & 1.45E-07 & 4.89E-10 & 4.07E-12 & 2.38E-12 \\ 
        \hline
    \end{tabular}
    \end{small}
\end{table}

\begin{figure}[htpb!]
\centering
\subfigure[SP]{
\label{Fig.sub.w1}
\begin{minipage}[t]{0.21\linewidth}
\centering
\includegraphics[width=1\textwidth]{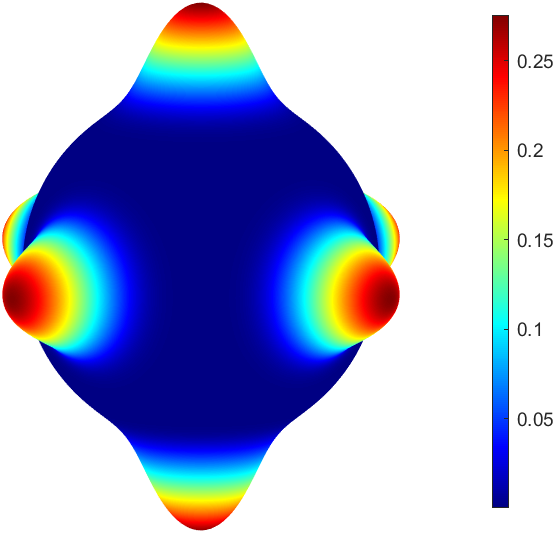}\vspace{1mm}
\quad
\includegraphics[width=0.95\textwidth]{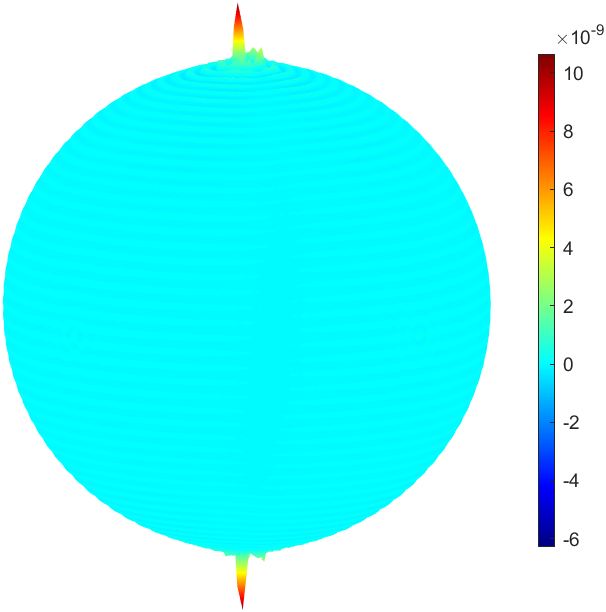}\vspace{1mm}
\quad
\includegraphics[width=1\textwidth]{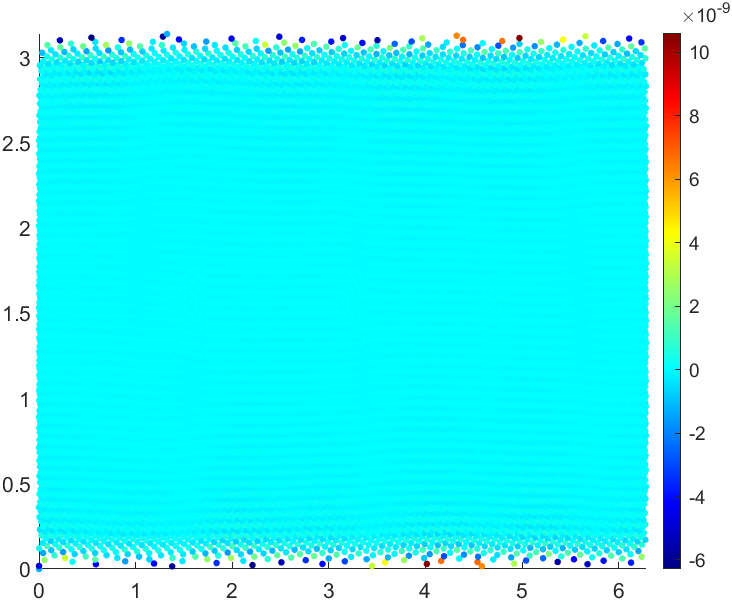}
\end{minipage}
}
\subfigure[SPD]{
\label{Fig.sub.w2}
\begin{minipage}[t]{0.21\linewidth}
\centering
\includegraphics[width=1\textwidth]{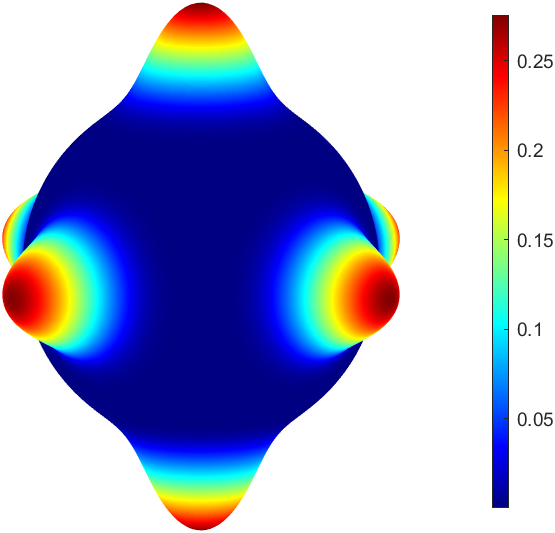}\vspace{1mm}
\quad
\includegraphics[width=1\textwidth]{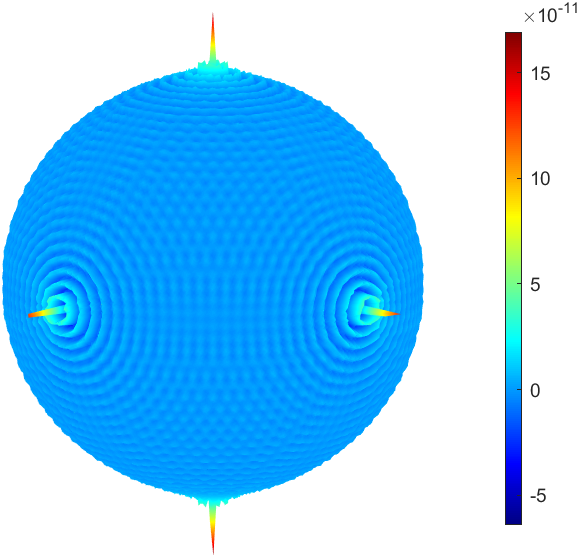}\vspace{1mm}
\quad
\includegraphics[width=1\textwidth]{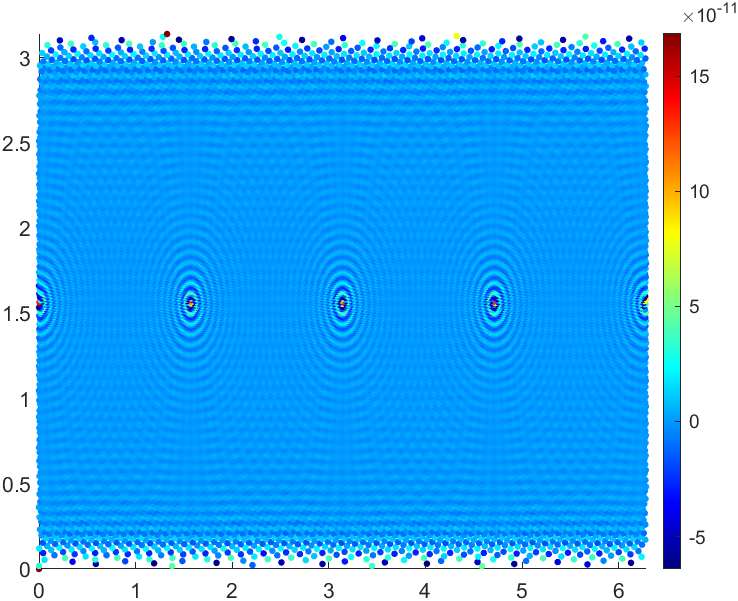}
\end{minipage}
}
\subfigure[UD]{
\label{Fig.sub.w3}
\begin{minipage}[t]{0.21\linewidth}
\centering
\includegraphics[width=1\textwidth]{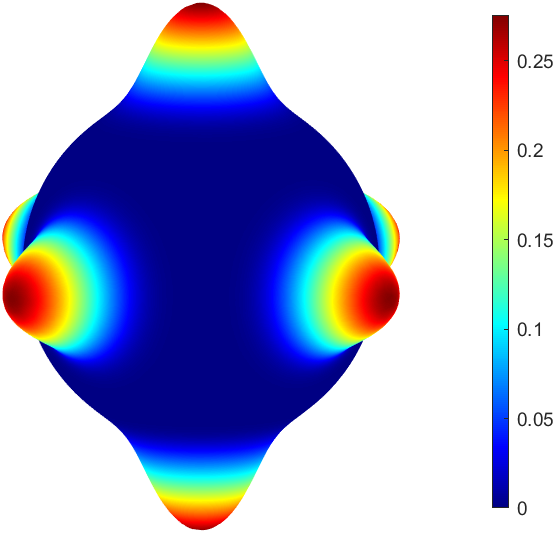}\vspace{3mm}
\quad
\includegraphics[width=0.95\textwidth]{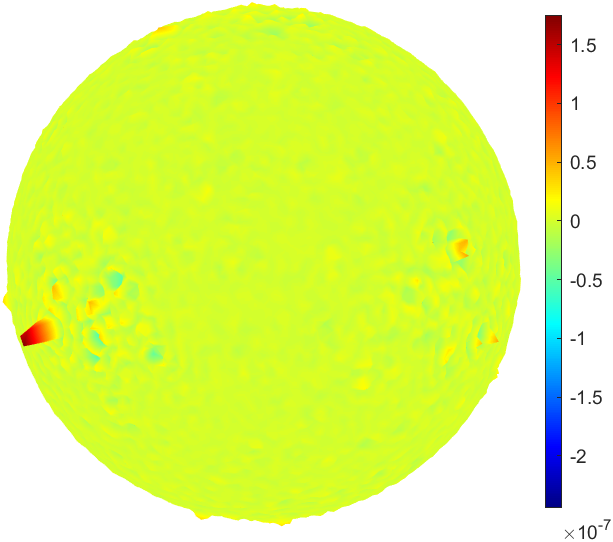}\vspace{3mm}
\quad
\includegraphics[width=1\textwidth]{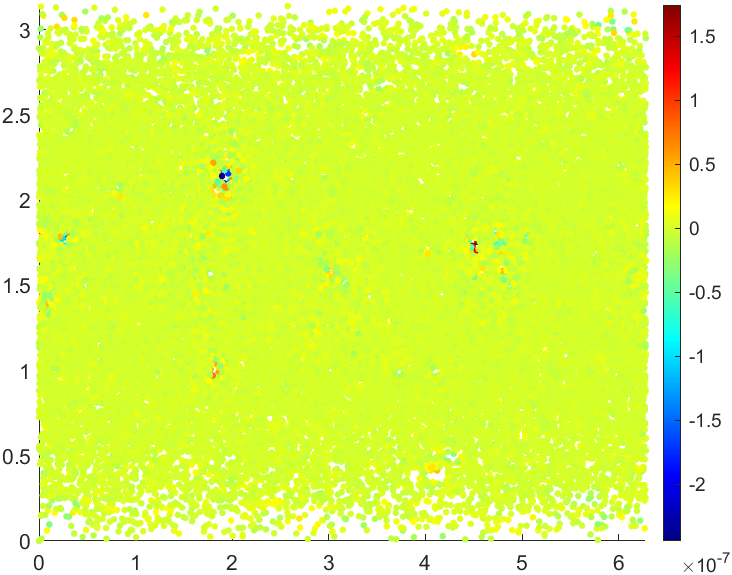}
\end{minipage}
}
\subfigure[SUD]{
\label{Fig.sub.w4}
\begin{minipage}[t]{0.21\linewidth}
\centering
\includegraphics[width=1\textwidth]{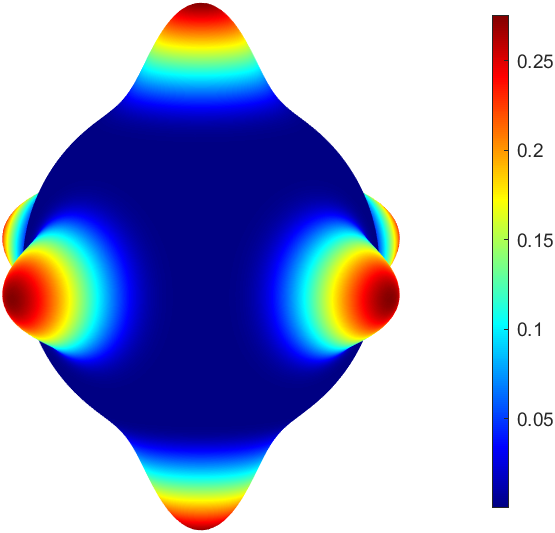}\vspace{3mm}
\quad
\includegraphics[width=0.95\textwidth]{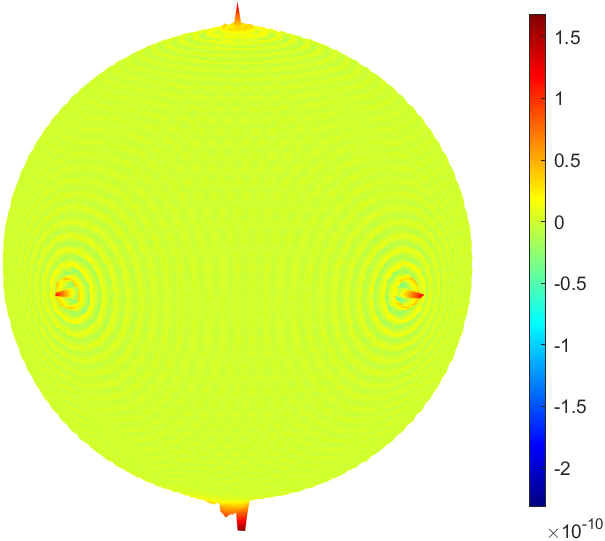}\vspace{3mm}
\quad
\includegraphics[width=1\textwidth]{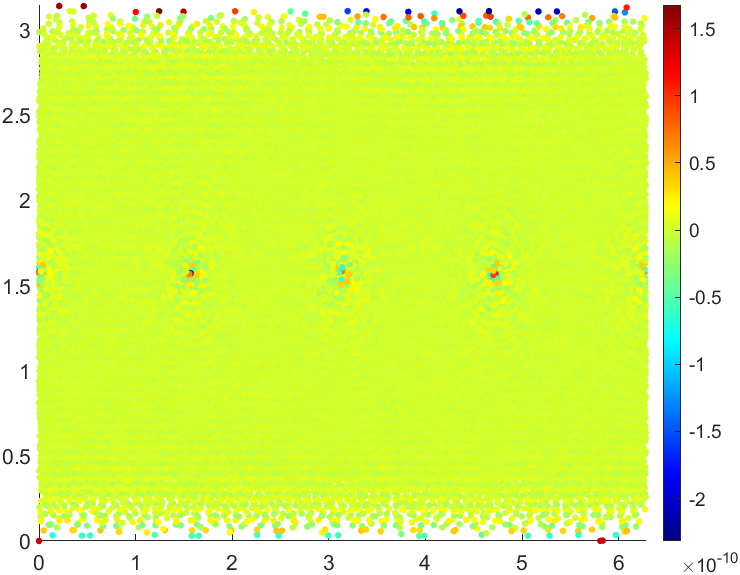}
\end{minipage}
}

\caption{Numerical simulation of real part of Projection term (first row), residual term (middle row), and the equirectangular projection of the residual (last row) for RBF $f_4$ under the setting $T=\frac{t}{2}$ and $W=\mathbf w$ on Algorithm~\ref{alg:projCG} for $t=200,N=(t+1)^2$.}
\label{Fig.labelw1}
\end{figure}

\begin{figure}[htpb!]
\centering
\subfigure[IV]{
\label{Fig.sub.w5}
\begin{minipage}[t]{0.21\linewidth}
\centering
\includegraphics[width=1\textwidth]{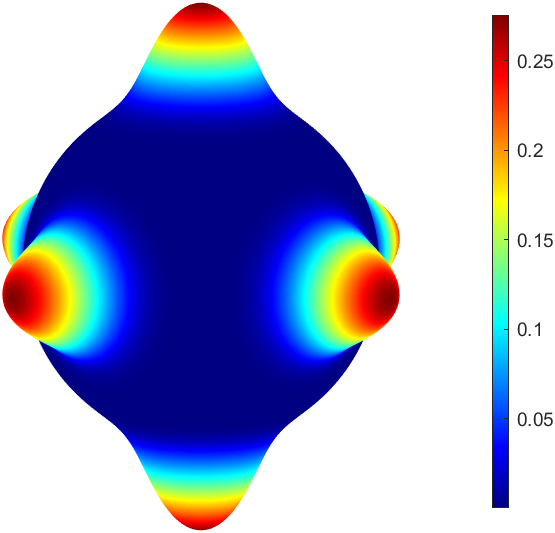}\vspace{2mm}
\quad
\includegraphics[width=0.95\textwidth]{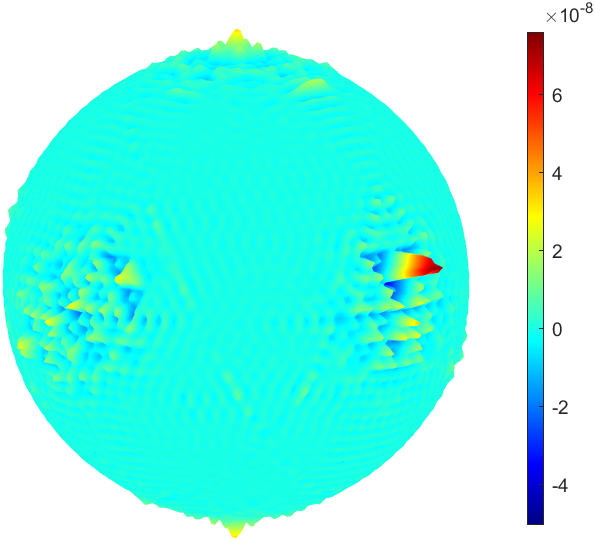}\vspace{2mm}
\quad
\includegraphics[width=1\textwidth]{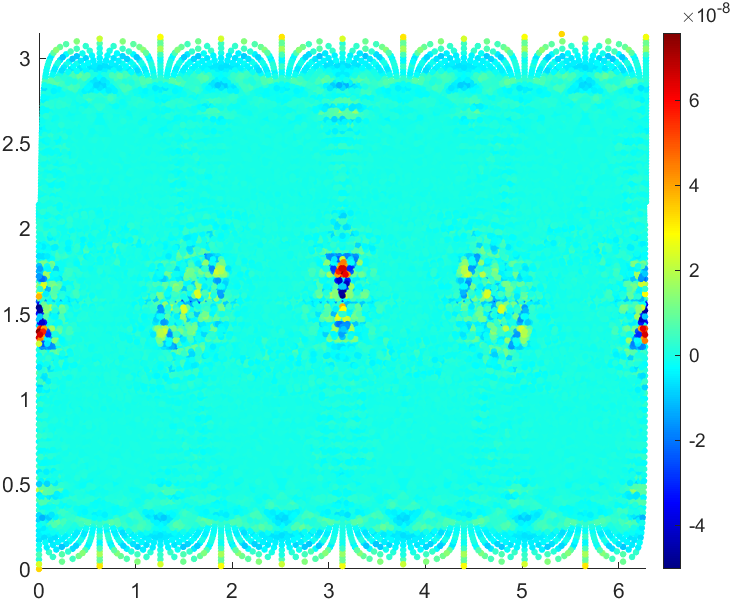}
\end{minipage}
}
\subfigure[SID]{
\label{Fig.sub.w6}
\begin{minipage}[t]{0.21\linewidth}
\centering
\includegraphics[width=1\textwidth]{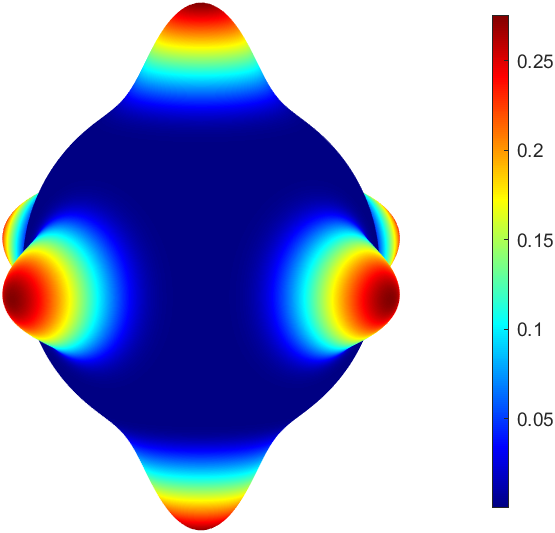}\vspace{1mm}
\quad
\includegraphics[width=1\textwidth]{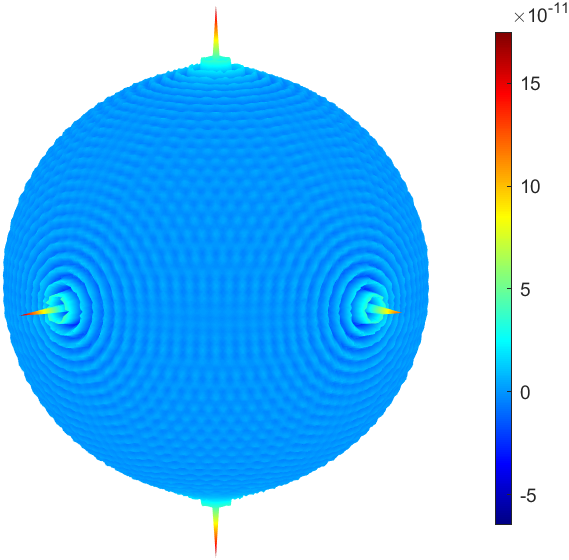}
\quad
\includegraphics[width=1\textwidth]{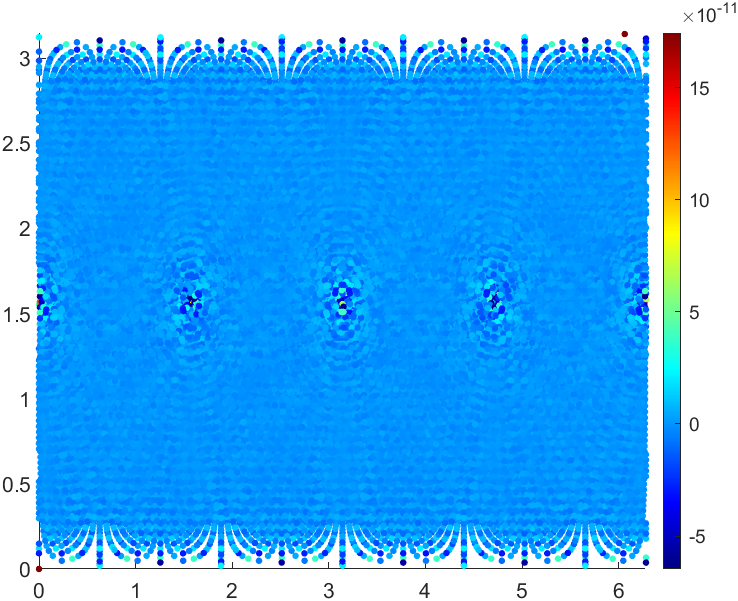}
\end{minipage}
}
\subfigure[HL]{
\label{Fig.sub.w7}
\begin{minipage}[t]{0.21\linewidth}
\centering
\includegraphics[width=1\textwidth]{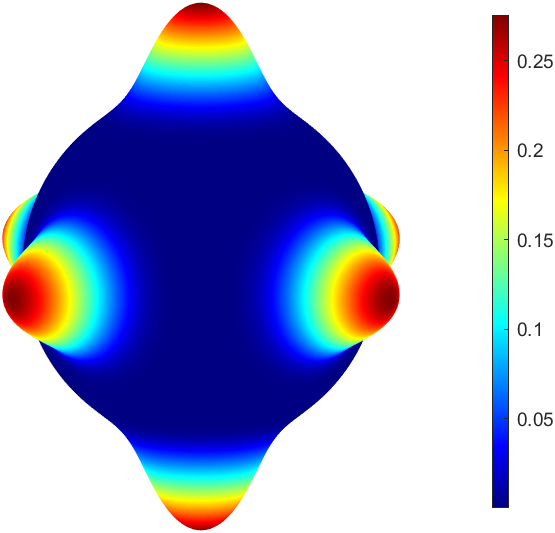}\vspace{3.5mm}
\quad
\includegraphics[width=0.94\textwidth]{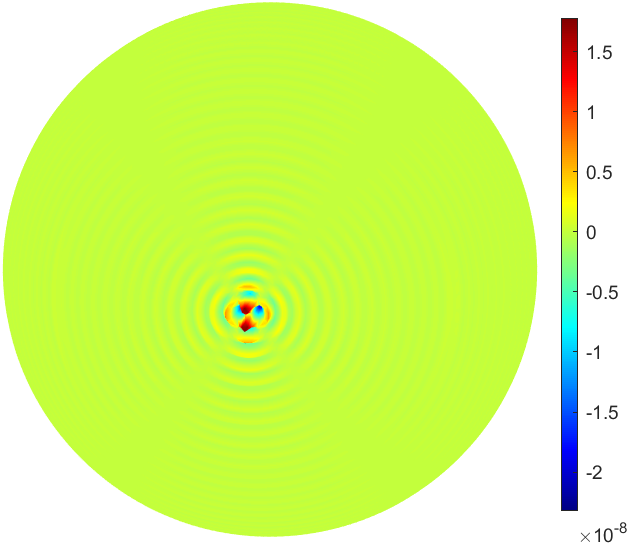}\vspace{3mm}
\quad
\includegraphics[width=1\textwidth]{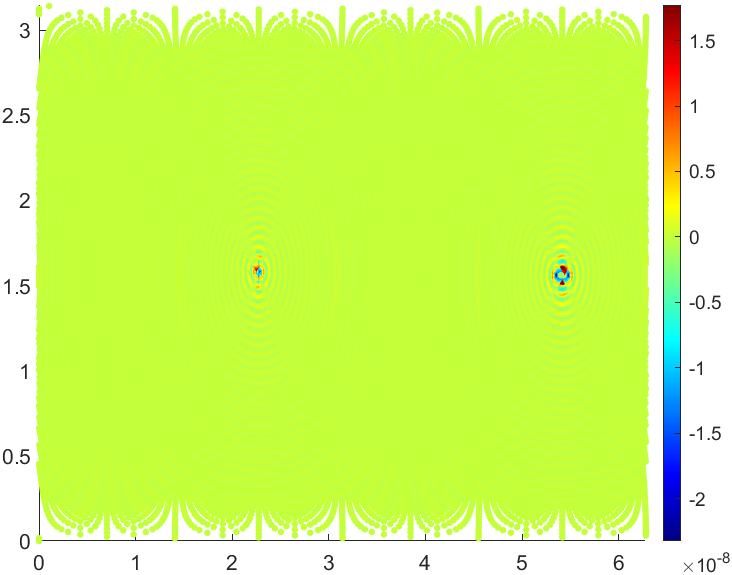}
\end{minipage}
}
\subfigure[SHD]{
\label{Fig.sub.w8}
\begin{minipage}[t]{0.21\linewidth}
\centering
\includegraphics[width=1\textwidth]{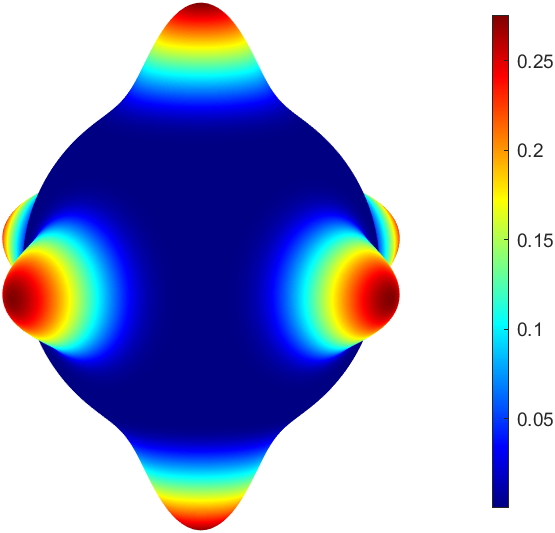}
\quad
\includegraphics[width=1\textwidth]{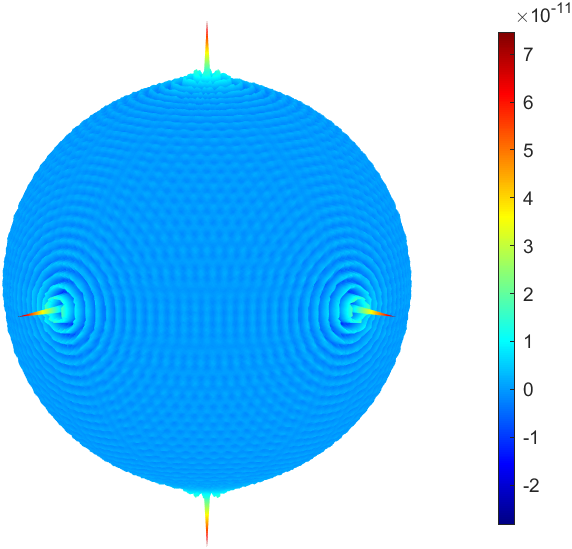}\vspace{1.5mm}
\quad
\includegraphics[width=1\textwidth]{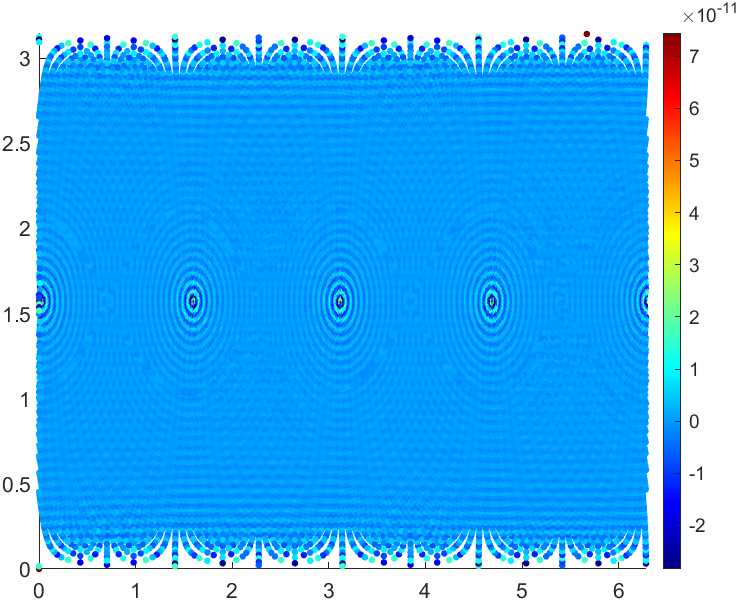}
\end{minipage}
}

\caption{Numerical simulation of the real part of Projection term (first row), the residual term (middle row), and the equirectangular projection of the residual (last row) for RBF $f_4$ under the setting $T=\frac{t}{2}$ and $W=\mathbf w$ on Algorithm~\ref{alg:projCG} for $t\approx 200,N\approx (t+1)^2$.}
\label{Fig.labelw2}
\end{figure}

From Tables~\ref{table:err2} and \ref{table:err4}, all point sets give good approximation error but the point sets with spherical $t$-design properties give a higher order of approximations. 
We display the related figures of the projection $v_1$ and residual $r$ of above point sets for RBF $f_4$ with degree $t\approx 200$ under the setting of Table~\ref{table:err2} on Figures~\ref{Fig.labelw1} and \ref{Fig.labelw2}. Comparing to related tables and figures, we can see that projection errors of the initial point sets under the setting of $T=\frac{t}{2}$ and $W=\sqrt{\mathbf w}$ are smaller than the setting of $T=\frac{t}{2}$ and $W=\mathbf w$, whereas the corresponding spherical $t$-designs point sets are almost the same (at about 1E-14 difference, but the results of $W=\sqrt{\mathbf w}$ is bit more smaller) except SUD at $t=400$ showed that when $W=\sqrt{\mathbf w}$ provides the smaller projection error.

%%%%%%%%%%%%%%%%%%%%%%%%%%%%%%%%%%%%%%
%%%%%%%%%%%%%%%%%%%%%%%%%%%%%%%%%%%%%%
\subsection{Denoising of Wendland functions by Projection}

\begin{figure}[htpb!]
%\centering
\subfigure[$\sigma=0.05$]{
\includegraphics[width=0.30\textwidth]{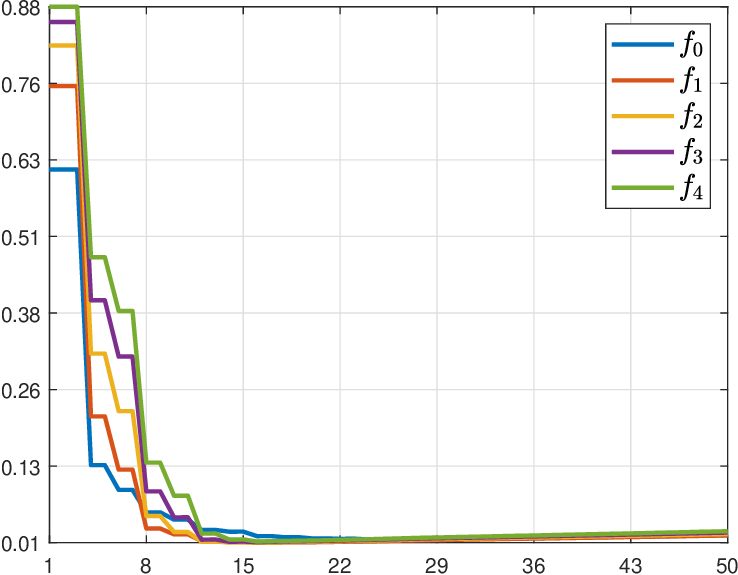}
}
\subfigure[$\sigma=0.1$]{
\includegraphics[width=0.30\textwidth]{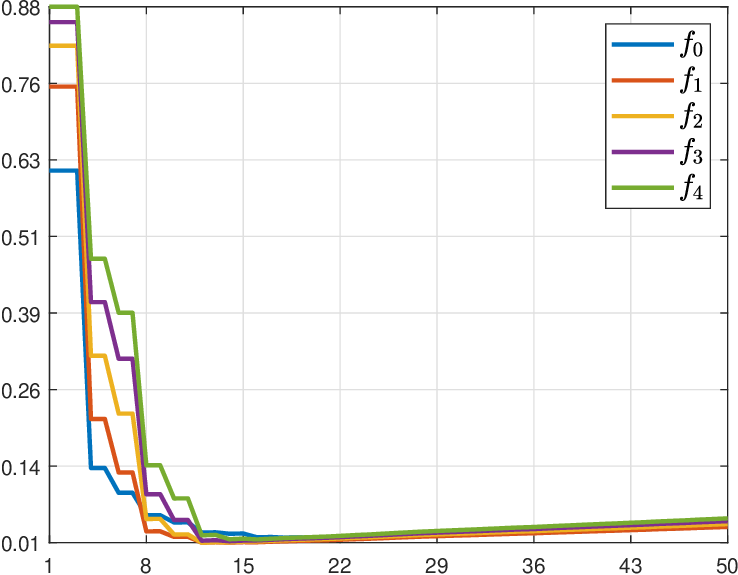}
}
\subfigure[$\sigma=0.15$]{
\includegraphics[width=0.30\textwidth]{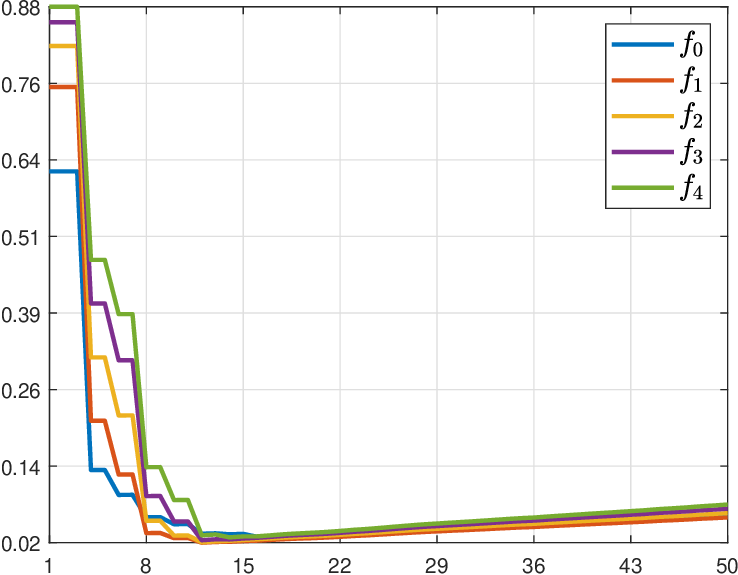}
}
\\
\subfigure[$\sigma=0.2$]{
\includegraphics[width=0.30\textwidth]{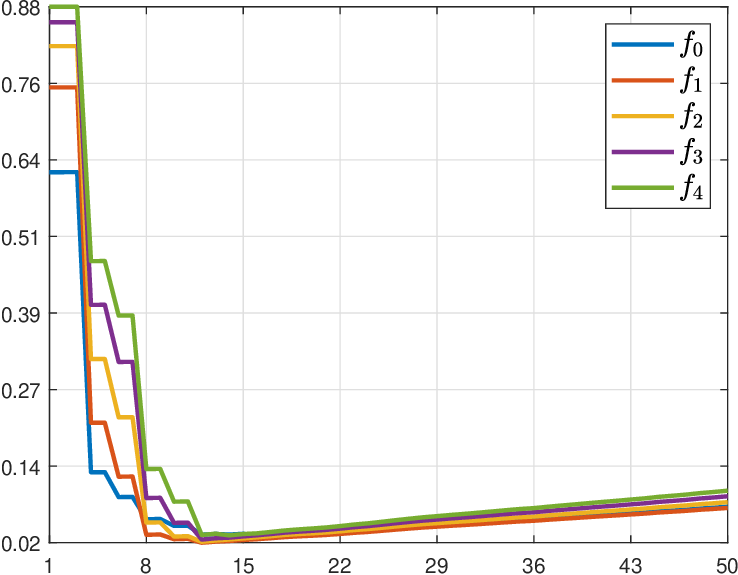}
}
\subfigure[$f_0$]{
\includegraphics[width=0.30\textwidth]{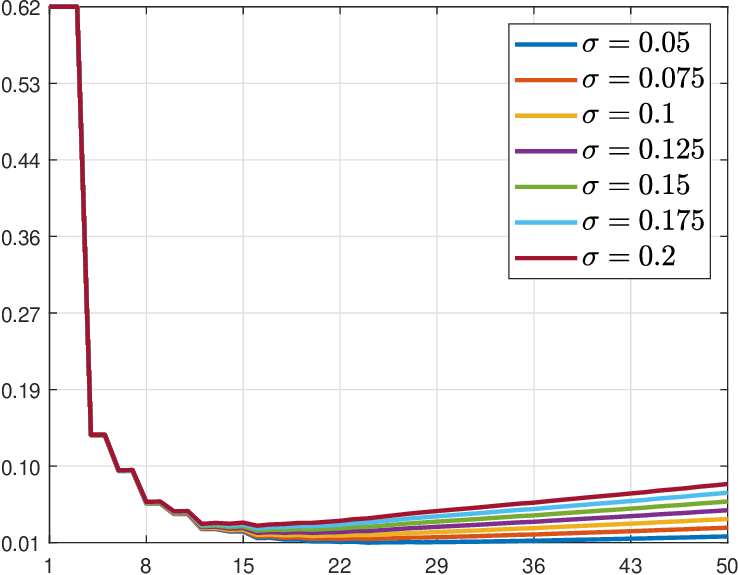}
}
\subfigure[$f_1$]{
\includegraphics[width=0.30\textwidth]{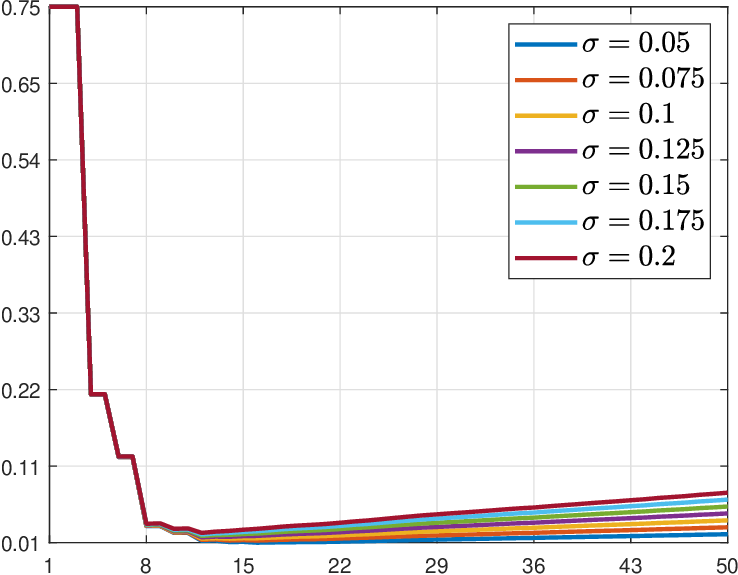}
}
\\
\subfigure[$f_2$]{
\includegraphics[width=0.30\textwidth]{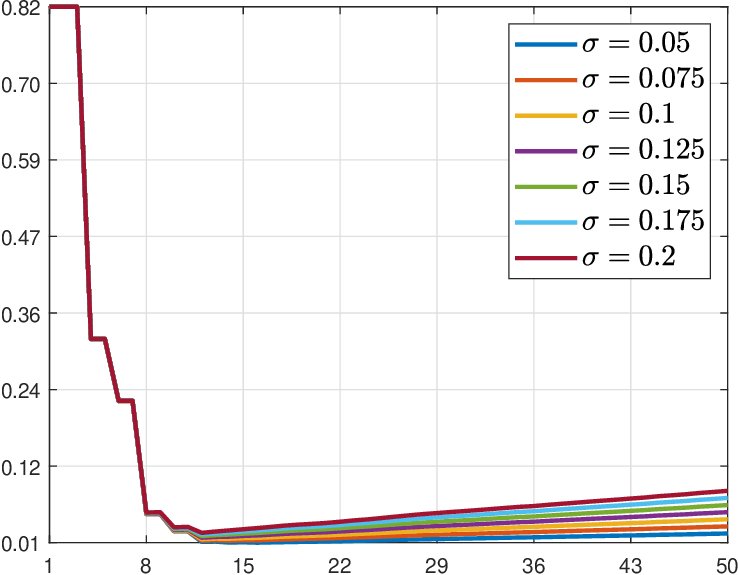}
}
\subfigure[$f_3$]{
\includegraphics[width=0.30\textwidth]{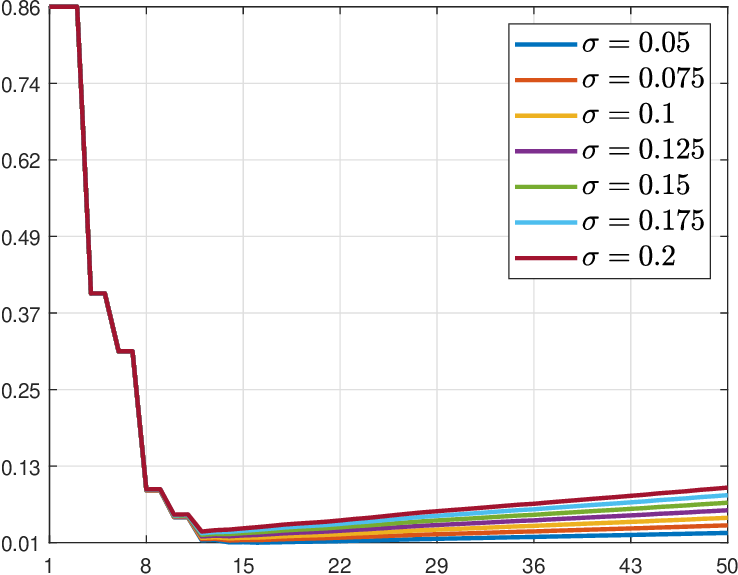}
}
\subfigure[$f_4$]{
\includegraphics[width=0.30\textwidth]{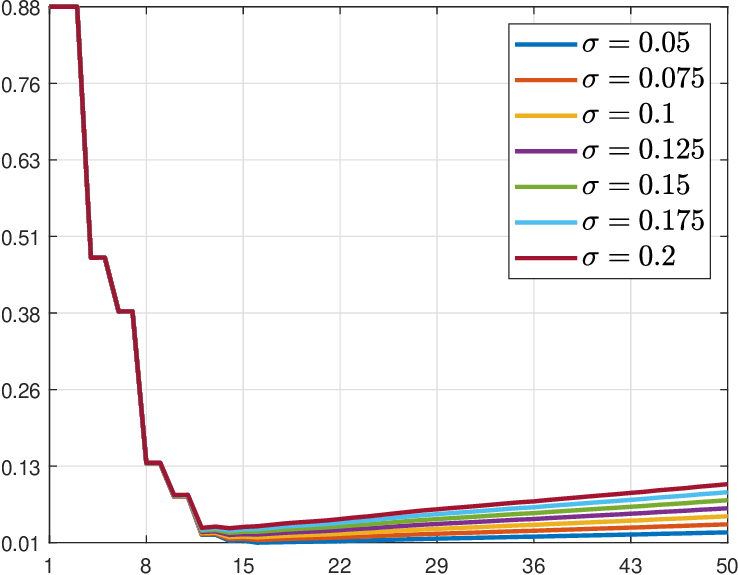}
}

\caption{Relative projection $L_2$-error with $err(f,f_k)$ for $k=0,1,2,3,4$ in different noise level $\sigma$. (a)-(d) are the extracted view from $\sigma=0.05,0.1,0.15,0.2$, $(e)-(i)$ are focusing on Wendland function $f_k$ for $k=0,1,2,3,4$.}
\label{Fig.bestpoly1}
\end{figure}

We then consider adding noise on Wendland function $f_k$ to observe the best polynomial degree for projecting noisy function on $\mathbb S^2$ since we have known the property of $f_k$. We generate Wendland function $f_k$ on the spherical $t$-design point set $X_N$ with $t=400,N=(t+1)^2$ from SPD. Let $f_{k,\sigma}(X_N)=f_k(X_N)+G_{\sigma}(X_N)$ be the noisy Wendland function generated by Gaussian white noise $G_{\sigma}$ with noise level $\sigma\lvert f_k\rvert_{\max}$. By applying Algorithm~\ref{alg:projCG} under the setting of maximum iterations $K_{\max}=1000$ and termination tolerance $\varepsilon$ = 2.2204e-16 (floating-point relative accuracy of MATLAB), and choosing projection degree $t$ from $1$ to $50$, we have $f_{k,\sigma}=f+g$, where $f\in\Pi_{t}$ and $g\notin\Pi_{t}$. Consider $\sigma$ taking from $0.05$ to $0.2$ with stepsize $0.025$, then we test the projection error of $f_k$ on $\mathbb S^2$ for $k=0,1,2,3,4$. We show relative projection errors for functions $f_k$ and $f$ on Figure~\ref{Fig.bestpoly1} and give the minimum error value for $f_k$ and $f$ on Table~\ref{table:bestpoly}.

We can see that the projection errors for $f_k$ and $f$ descend rapidly and reach the minimum value at certain degrees, then climb up slowly compared with the rate of descent. The projection errors for $f_{k,\sigma}$ and $f$ also descend rapidly for the same degree but will converge to 0 slowly as the degree is sufficiently large. From Table~\ref{table:bestpoly}, we obtain the best projection degrees for projecting $f_k$ from noisy function $f_{k,\sigma}$, which mostly fall on 12 to 16. Thus, we can choose the appropriate spherical $t$-design point sets with degree $t$ at about 12 to 16.
% as final level for semi-discrete spherical framelet systems to denoise in the next section.

\begin{table}[htpb!]
\centering
\caption{The minimum relative projection $L_2$-error $err(f,f_k)$  along with corresponding degree $t$ (inside the parentheses) in different noise level $\sigma$.}\label{table:bestpoly}
\begin{small}
\begin{tabular}{llllll}
\hline
$\sigma$ & $f_0$ & $f_1$ & $f_2$ & $f_3$ & $f_4$ \\ \hline
0.05 & 1.28E-02(24) & 7.40E-03(16) & 7.60E-03(16) & 8.40E-03(16) & 9.31E-03(16) \\ \hline
0.075 & 1.68E-02(24) & 1.05E-02(16) & 1.07E-02(14) & 1.20E-02(14) & 1.37E-02(16) \\ \hline
0.1 & 2.03E-02(20) & 1.30E-02(14) & 1.28E-02(12) & 1.53E-02(14) & 1.82E-02(16) \\ \hline
0.125 & 2.38E-02(20) & 1.49E-02(12) & 1.51E-02(12) & 1.84E-02(12) & 2.19E-02(14) \\ \hline
0.15 & 2.64E-02(16) & 1.69E-02(12) & 1.75E-02(12) & 2.08E-02(12) & 2.55E-02(14) \\ \hline
0.175 & 2.91E-02(16) & 1.89E-02(12) & 2.00E-02(12) & 2.34E-02(12) & 2.92E-02(14) \\ \hline
0.2 & 3.19E-02(16) & 2.11E-02(12) & 2.25E-02(12) & 2.61E-02(12) & 3.29E-02(14) \\ \hline
\end{tabular}
\end{small}
\end{table}

%%%%%%%%%%%%%%%%%%%%%%%%%%%%%%%%%%%%%%
%%%%%%%%%%%%%%%%%%%%%%%%%%%%%%%%%%%%%%
\subsection{Approximation of non-smooth functions}
\begin{figure}[htpb!]
\centering
\subfigure[IV]{
\label{Fig.sub.ns3}
\includegraphics[width=0.33\textwidth]{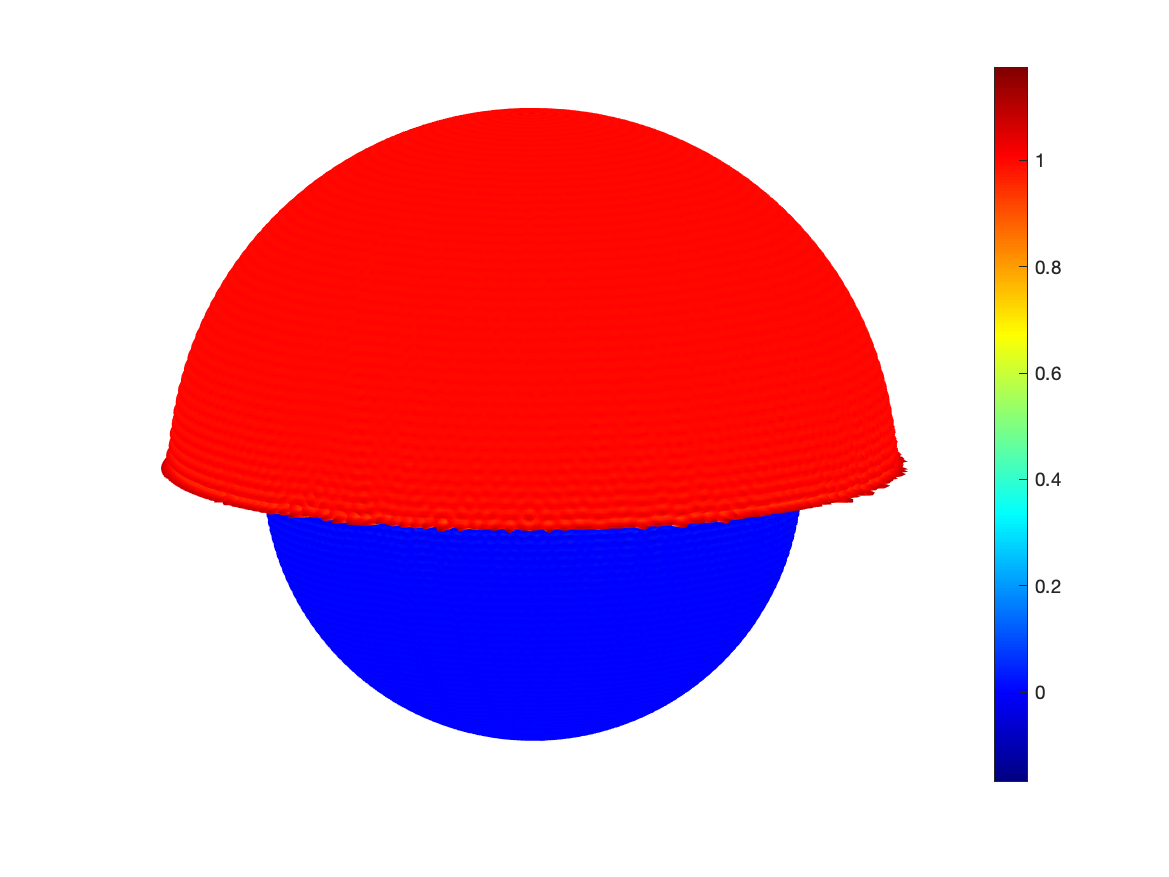}
\includegraphics[width=0.33\textwidth]{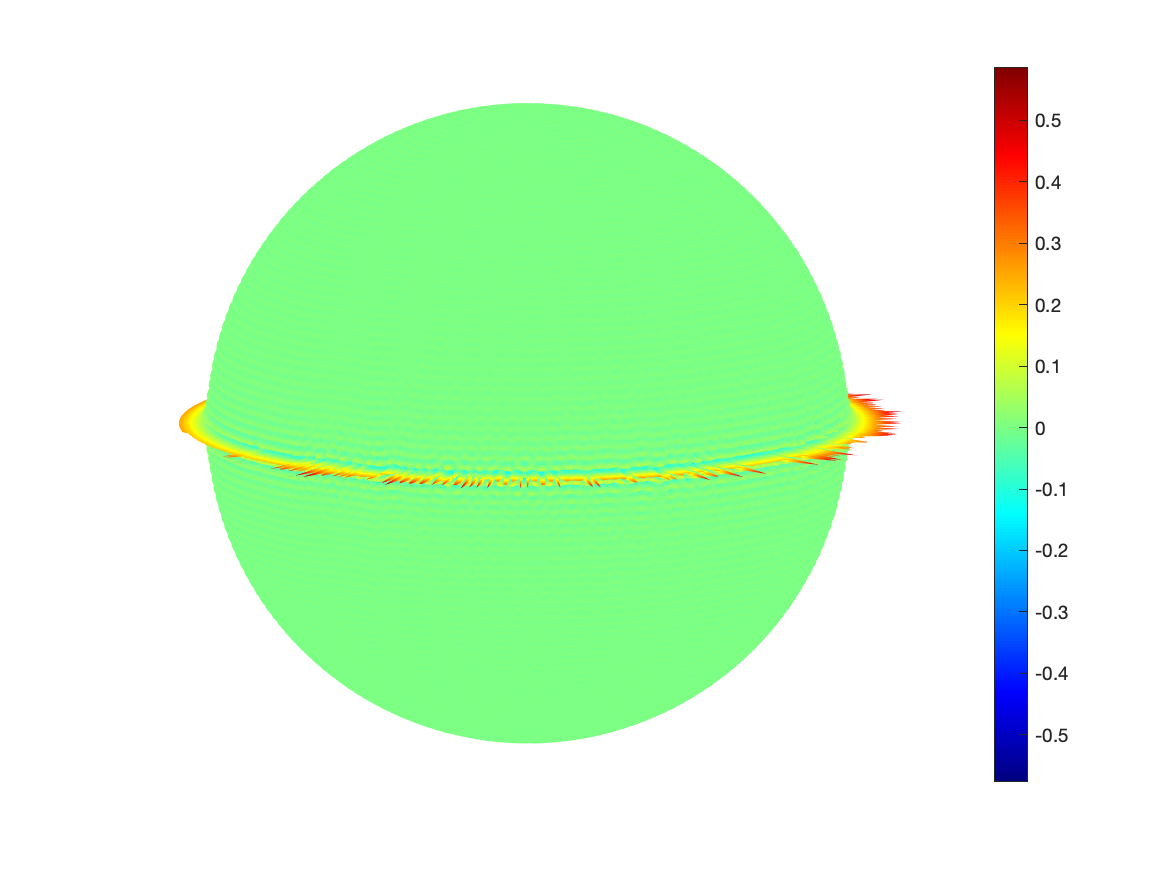}
\includegraphics[width=0.33\textwidth]{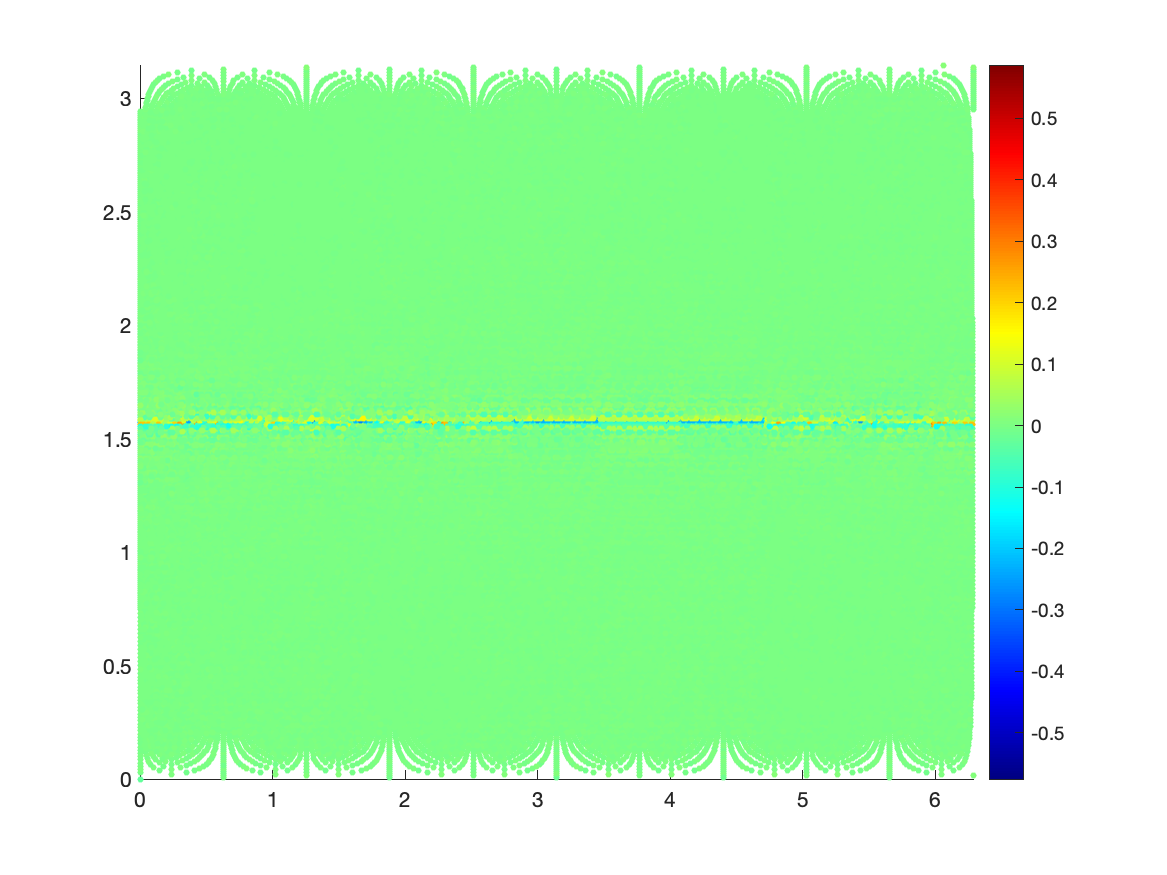}
}
\subfigure[SID]{
\label{Fig.sub.ns4}
\includegraphics[width=0.33\textwidth]{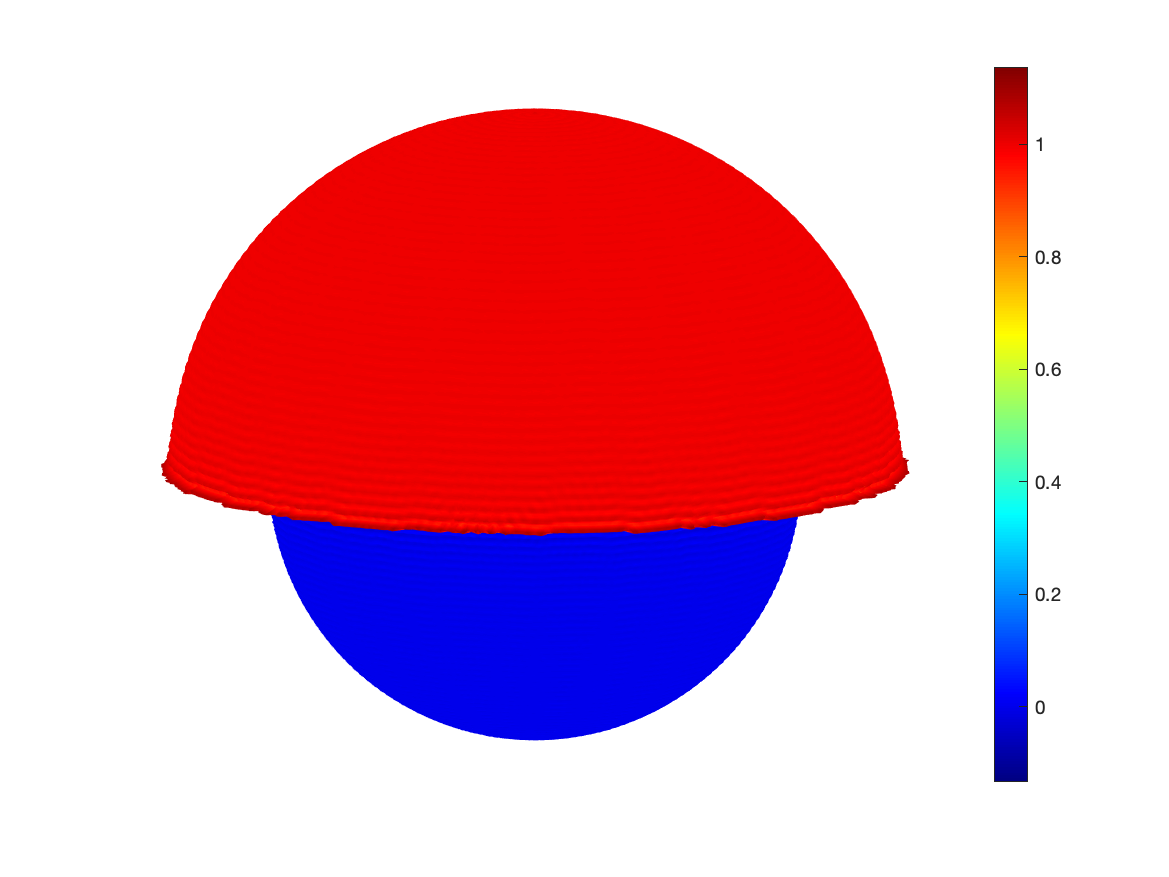}
\includegraphics[width=0.33\textwidth]{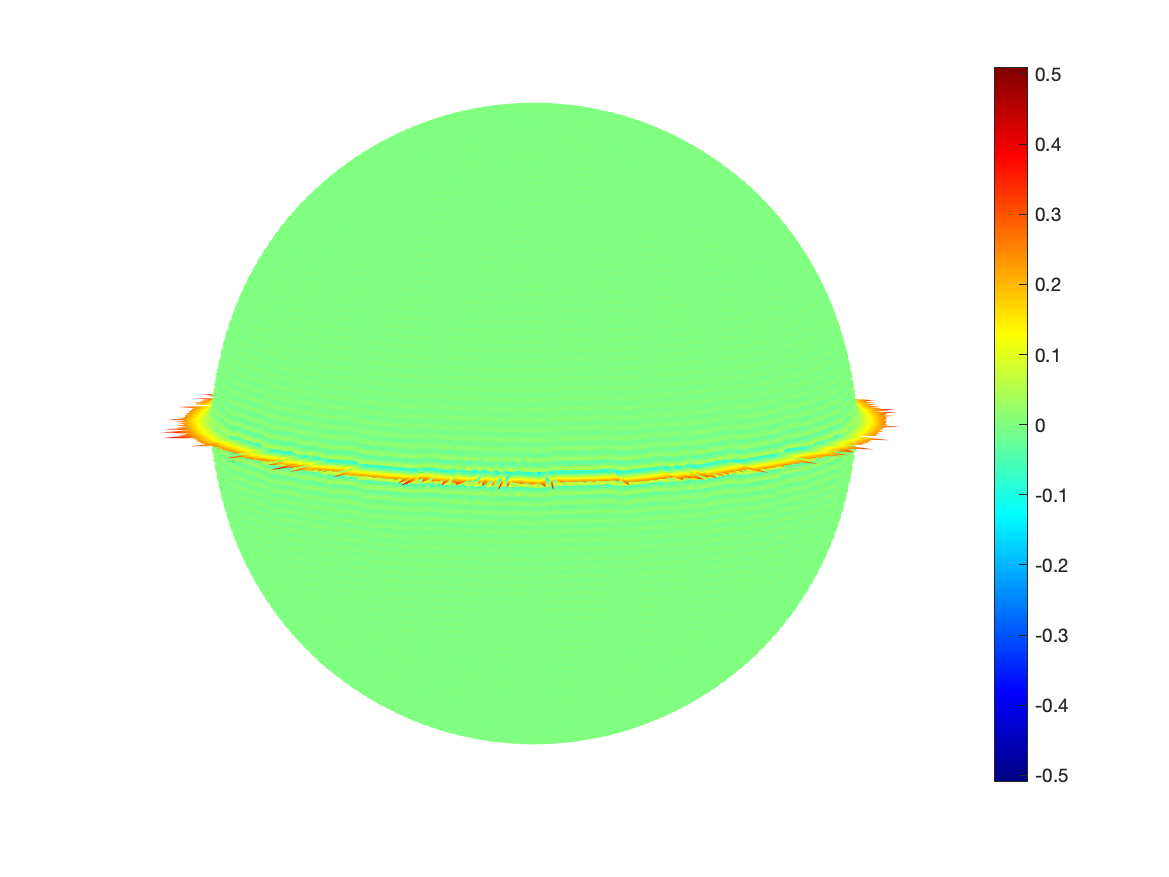}
\includegraphics[width=0.33\textwidth]{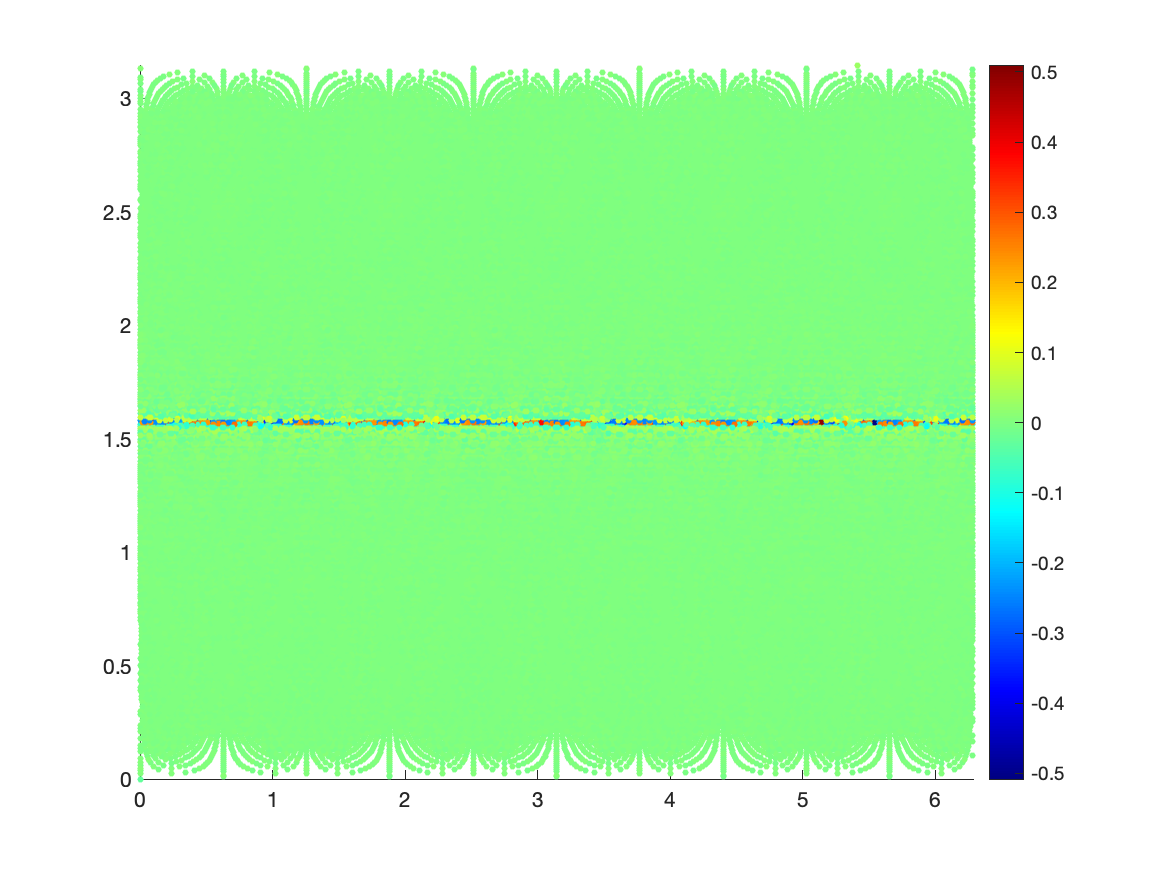}
}

\caption{Numerical simulation of real part of Projection term (left column), residual term (middle column), and the equirectangular projection of the residual (right column) for spherical $1_{\mathcal D}$ function under the setting $T=\frac{t}{2}$ and $W=\sqrt{\mathbf w}$ on Algorithm~\ref{alg:projCG} for $t=403, N=163842$.}
\label{Fig.label.ns2}
\end{figure}

Finally,  instead of considering differentiable functions in Section~\ref{subsec:smth}, we consider here non-smooth functions approximation with spherical designs. We define a \emph{Spherical $\bm 1_{\mathcal D}$ function}  as follows
\begin{equation*}\label{1dfcn}
f_{\bm 1_{\mathcal D}}:=\begin{cases}
  1, & \bm x\in\mathcal D,\\
  0, & \bm x\notin\mathcal D,
  \end{cases}
\end{equation*}
where $\mathcal D:=\{(\theta,\phi)\setsep\theta\in [0,\frac{\pi}{2}],\phi\in [0,2\pi);(\theta,\phi)\in\mathbb S^2 \}$ is the northern hemisphere. Similarly, we set $T=\frac{t}{2}$ and weight $W=\mathbf w$, and $T=\frac{t}{2}$ and weight $W=\sqrt{\mathbf w}$ on Algorithm~\ref{alg:projCG} to obtain Table~\ref{table:err1}. We can see from Table~\ref{table:err1}, that the higher degree $t$, the lower the error of the projection of spherical $1_{\mathcal D}$ function. Besides, there is no significant difference between the two choices of weights. In addition, only a little refinement of the error for point set SPD with degree $t=400$ and SID with $t=201,403$ than their corresponding initial point sets. We also display the example of projection in Figure~\ref{Fig.label.ns2}. One can see that the discontinuity of the function $f_{\bm 1_{\mathcal D}}$ along the equator brings the main challenge for the polynomial approximation. 

\begin{table}[htpb!]
    \centering
    \caption{Relative $L_{2}$-errors $err(f_{w},f_{\bm 1_{\mathcal D}})$ and $err(f_{\sqrt w},f_{\bm 1_{\mathcal D}})$ for projecting spherical $1_{\mathcal D}$ function. $f_{w}$ and $f_{\sqrt w}$ represent the relative error with functions on point sets along with corresponding degree $T=\frac{t}{2}$ take weights $W=\mathbf w$ and $W=\sqrt{\mathbf w}$, respectively on Algorithm~\ref{alg:projCG}.}\label{table:err1}
    \begin{small}
    \begin{tabular}{lllll}
    \hline
        $t$ & $N$ & $Q_N$ & $f_{w}$ & $f_{\sqrt w}$ \\
        \hline
        200 & 40401 & SP & 5.46E-02 & 5.46E-02 \\
        200 & 40401 & SPD & 5.46E-02 & 5.46E-02 \\
        400 & 160801 & SP & 3.89E-02 & 3.89E-02 \\ 
        400 & 160801 & SPD & 3.88E-02 & 3.88E-02 \\ 
        \hline
        200 & 40401 & UD & 4.85E-02 & 4.85E-02 \\
        200 & 40401 & SUD & 5.21E-02 & 5.21E-02 \\
        400 & 160801 & UD & 3.52E-02 & 3.52E-02 \\ 
        400 & 160801 & SUD & 3.68E-02 & 3.68E-02 \\ 
        \hline
        201 & 40962 & IV & 5.50E-02 & 5.50E-02 \\
        201 & 40962 & SID & 5.40E-02 & 5.40E-02 \\
        403 & 163842 & IV & 3.84E-02 & 3.84E-02 \\ 
        403 & 163842 & SID & 3.65E-02 & 3.65E-02 \\ 
        \hline
        220 & 49152 & HL & 5.07E-02 & 5.07E-02 \\
        220 & 49152 & SHD & 5.07E-02 & 5.07E-02 \\
        442 & 196608 & HL & 3.56E-02 & 3.56E-02 \\ 
        442 & 196608 & SHD & 3.56E-02 & 3.56E-02 \\ 
        \hline
    \end{tabular}
    \end{small}
\end{table}

%%%%%%%%%%%%%%%%%%%%%%%%%%%%%%%%%%%%%%
%%%%%%%%%%%%%%%%%%%%%%%%%%%%%%%%%%%%%%
%%%%%%%%%%%%%%%%%%%%%%%%%%%%%%%%%%%%%%
\section{Spherical framelets for spherical signal processing using spherical designs}
\label{sec:apply}
In this section, we discuss the further application of spherical designs in spherical framelet constructions and spherical signal processing. In general, given a group of spherical design point sets, for a noisy function $f_{\sigma}$ with  Gaussian white noise of noise level $\sigma\lvert f\rvert_{\max}$,  following the orthogonal projection procedure of Section~\ref{subsec:proj}, we have $f_{\sigma}=f+g$, where the projected function $f\in\Pi_{t_{J}}$ and residual $g\notin\Pi_{t_{J}}$. Notice that both projected function $f$ and residual $g$ have noise.  Then we use spherical $t$-design point sets to decompose $f$ by semi-discrete truncated spherical tight framelet system $\mathcal{F}_{J_0}^J(\eta,\mathcal Q)$. After that, we use threshold methods based on spherical caps to denoise the framelet coefficients of $f$, and then we reconstruct the function by the denoised coefficients to obtain $f_{thr}$. Similarly, $g$ is also being denosed to obtain $g_{thr}$. Finally, we get the denoised function $F_{thr}=f_{thr}+g_{thr}$. Below, we introduce the truncated spherical framelet systems, the spherical caps, the thresholding techniques, and the denoising experiments on Wendland functions.

\subsection{Semi-discrete spherical framelet systems}
\label{subsec:fmt}
First, we briefly introduce the semi-discrete spherical framelet systems based on spherical harmonics and spherical design point sets. 

We start with a filter bank
$\eta = \{a; b_1,\ldots, b_n\}\subset l_1(\mathbb{Z}):=\{h=\{ h_k\}_{k\in\mathbb Z}\subset\mathbb C\setsep\sum_{k\in\mathbb Z}\lvert h_k\rvert<\infty\}$ satisfying the partition of unity condition:
\begin{align}
\label{PUC:eta}
|\hat a(\xi)|^2+\sum_{s\in[n]} |\hat b_s(\xi)|^2=1,\quad \xi\in\R.
\end{align}
where for a filter (mask) $h=\{h_k\}_{k\in\mathbb Z}\subset\mathbb C$, its Fourier series  $\hat h$ is given by $\hat h(\xi):=\sum_{k\in\mathbb Z}h_k \mathrm e^{-2\pi\mathrm ik\xi}$, for $\xi\in\mathbb R$,  and for a positive integer $n$,  we denote $[n]:=\{1,\ldots, n\}$.

The \emph{quadrature (cubature) rule} $Q_{N_j}=(X_{N_j},\mathbf w_{j})$ on $\mathbb S^2$ at scale $j$  is a collection of point set $X_{N_j}:=\{\bm x_{j,k}\setsep k\in[N_j]\}$ and weight $\mathbf w_{j}:=\{w_{j,k}\setsep k\in[N_j]\}$, where $N_j$ is the number of points at scale $j$. A quadrature rule $Q_{N_j}$ is \emph{polynomial-exact} up to degree $t_j\in\mathbb N_0$ if $\sum_{k=1}^{N_j} w_{j,k} p(\bm x_{j,k})=\int_{\mathbb S^2}p(\bm x)\mathrm d\mu_2(\bm x)$  for all $p\in\Pi_{t_j}$. We denote a polynomial-exact quadrature rule of degree $t_j$ as $Q_{N_j}=:Q_{N_j,t_j}$. The spherical $t$-design  $X_N=\{\bm x_1,\ldots, \bm x_N\}$ forms a polynomial-exact quadrature rule $Q_{N,t}:=(X_N, \mathbf w)$ of degree $t$ with weight $\mathbf w\equiv \frac{4\pi}{N}$.
We consider
$\mathcal{Q}:=\mathcal{Q}_{J_0}^{J+1}:=\{Q_{N_j,t_j}:j=J_0,\ldots,J+1\}$ be a set of polynomial-exact quadrature rules satisfying  $t_{j+1}=2t_j$.

For a function $f\in L^1(\R)$, its Fourier transform $\hat f$ is defined by $\hat f(\xi):=\int_\R f(x)e^{-2\pi \mathrm i x\xi}dx$. For a fixed fine scale $J\in\Z$, we set
\begin{align}
\label{alphaJ1}
\hat\alpha^{(J+1)}(\frac{\ell}{t_{J+1}})
=
\begin{cases}
 1  & \mbox{for } \ell\le t_J, \\
 \\
 0  & \mbox{for } \ell>t_J,
\end{cases}
\end{align}
and we recursively define $\hat \alpha^{(j)},\hat \beta_s^{(j)}$ from $\hat\alpha^{(j+1)}$ by
\begin{align}
\label{eta2Psi1}
\hat \alpha^{(j)}(\frac{\ell}{t_j})
=&\hat \alpha^{(j)}(2\frac{\ell}{t_{j+1}})
=\hat a(\frac{\ell}{t_{j+1}})\hat \alpha^{(j+1)}(\frac{\ell}{t_{j+1}}), \\
\label{eta2Psi2}
\hat \beta_s^{(j)}(\frac{\ell}{t_j})
=&\hat \beta_s^{(j)}(2\frac{\ell}{t_{j+1}})
=\hat b_s(\frac{\ell}{t_{j+1}})\hat \alpha^{(j+1)}(\frac{\ell}{t_{j+1}}),\quad s\in[n],
\end{align}
for $j$ decreasing from $J$ to $J_0$.  Then, we obtain
\begin{align}
\label{Psi:new}
\Psi=\{\alpha^{(j)},\beta_s^{(j)}\setsep j=J_0,\ldots,J; s\in[n]\}.
\end{align}

%For convenience, we define the index set
%\begin{align}
%\label{def:indexSet}
%\mathcal I_t:=\{(\ell,m)\setsep \ell=0,\ldots,t; m=-\ell,\ldots,\ell\}.
%\end{align}
Using the above $\Psi$ and $\mathcal{Q}$, we can define the {\em truncated (semi-discrete) spherical framelet system} $\mathcal{F}_{J_0}^{J}(\eta,\mathcal{Q})$  from the spherical designs as
\begin{align}
\label{def:sphTruncated}
\mathcal{F}_{J_0}^{J}(\eta,\mathcal Q):=\{\varphi_{J_0,k}\setsep k\in [N_{J_0}]\}\cup\{\psi_{j,k}^{(s)}\setsep k\in [N_{j+1}],s\in[n]\}_{j=J_0}^{J},
\end{align}
where  $\varphi_{j,k}$ and $\psi_{j,k}^{(s)}$ are defined by
\begin{align}
\label{eq24.1}
\varphi_{j,k}(\bm x)&:=\sqrt{w_{j}}\sum_{(\ell,m)\in\mathcal{I}_{t_j}} \hat\alpha^{(j)}(\frac{\ell}{t_j})\overline{Y_{\ell}^{m}(\bm x_{j,k})}Y_{\ell}^{m}(\bm x),\\
\label{eq25.1}
\psi_{j,k}^{(s)}(\bm x)&:=\sqrt{w_{j+1}}\sum_{(\ell,m)\in\mathcal{I}_{t_{j+1}}} \hat{\beta}_{s}^{(j)}(\frac{\ell}{t_j})\overline{Y_{\ell}^{m}(\bm x_{j+1,k})}Y_{\ell}^{m}(\bm x).
\end{align}
The truncated system $\mathcal{F}_{J_0}^J(\eta,\mathcal{Q})$ is completely determined by the filter bank $\eta$ and the quadrature rules $\mathcal{Q}$.
It was shown in \cite[Theorem 4.3]{xiao2023spherical} that $\mathcal{F}_{J_0}^J(\eta,\mathcal{Q})$  is a tight frame for $\Pi_{t_J}$. That is, for all $f\in\Pi_{t_J}$, we have 
\[
f = \sum_{k=1}^{N_{0}}\mathpzc v_{J_0,k}\varphi_{J_0,k}+\sum_{j=J_0}^J\sum_{k=1}^{N_{j}}\sum_{s=1}^n \mathpzc w_{j,k}^{(s)}\psi_{j,k}^{(s)},
\] 
where $\mathpzc v_{j,k}:=\langle f,\varphi_{j,k}\rangle_{L_2(\mathbb{S}^2)}$ and  $\mathpzc w_{j,k}^{(s)}:=\langle f,\psi_{j,k}^{(s)}\rangle_{L_2(\mathbb{S}^2)}$. 

Therefore, we can use the spherical tight framelet system to decompose $f$ into framelet coefficient sequences  $\{\mathpzc v_{j,k};\mathpzc w_{j,k}^{(s)}:s=1,\ldots,n\}$. We then apply the thresholding techniques for denoising the framelet coefficient sequence.

\subsection{Spherical caps}
\label{subsec:caps}
\begin{figure}[htpb!]
\centering

\subfigure[Partial view]{
\includegraphics[width=0.46\textwidth]{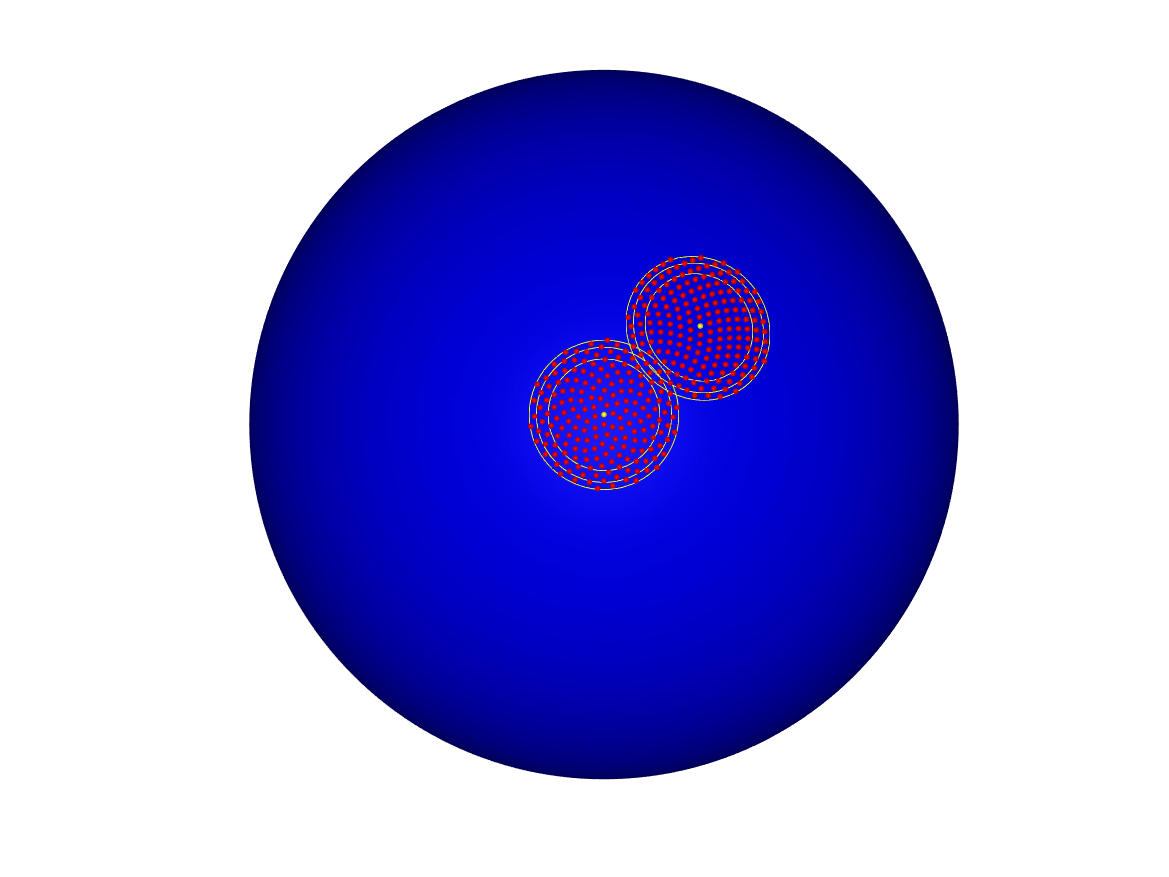}
}
\quad
\subfigure[General view]{
\includegraphics[width=0.46\textwidth]{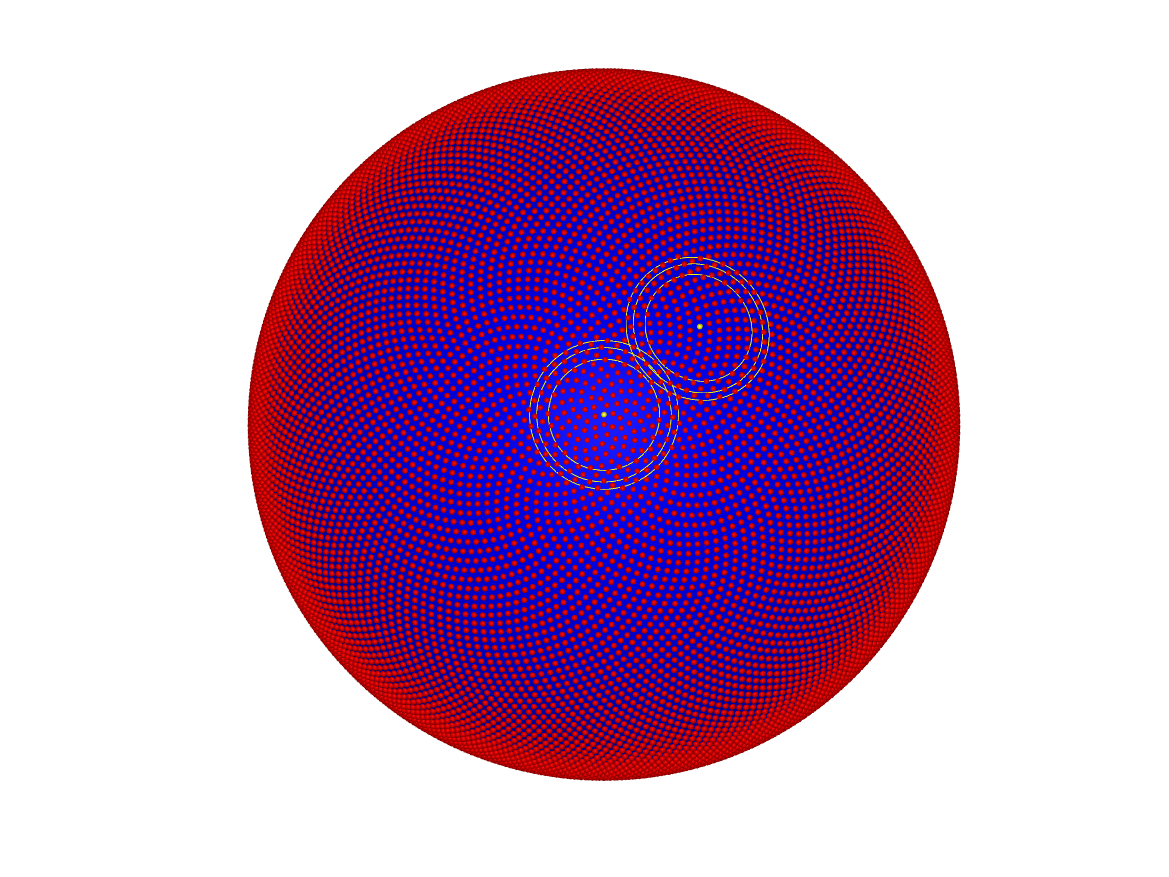}
}

\caption{Spherical caps (rnn-based) on spherical $128$-design point set (SPD) with $N=(t+1)^2$ for cap layer $i=15,22,27$ on the centroid $\bm x_1$ and $\bm x_{600}$ points}
\label{Fig.cap}
\end{figure}

Given a point set $X_N\subset\mathbb S^2$, for each point $\bm x\in X_N$, We want to get geodesic information for the connection of points on $X_{N_j}$. Naturally, we consider the spherical cap. In one way, we can generate index sets of $\bm x$ based on \emph{radius nearest neighbor (rnn)} $\mathcal C_{R}$. To get this, we firstly define a spherical cap $C(\bm x,r)\subset\mathbb S^2$ centred at $\bm x\in\mathbb S^2$ with radius $r\in\mathbb N_0$ is given by the set
\begin{align*}
C(\bm x,r):=\{\bm y\in\mathbb S^2:\lVert\bm x\times\bm y\rVert\leq r\},
\end{align*}
where $\lVert\bm x\times\bm y\rVert$ is the magnitude of the cross product of $\bm x$ and $\bm y$. Then the boundary of $C(\bm x,r)$ is
\begin{align*}
\partial C(\bm x,r):=\{\bm y\in\mathbb S^2:\lVert\bm x\times\bm y\rVert= r\}.
\end{align*}
Then we have $\{\bm x_1,\ldots,\bm x_s\}\in C(\bm x,r)$ where $s$ is the total number of points found in $C(\bm x,r)$. Thus, we obtain index set $\mathcal C_{R}(X_N)$ of points in $C(\bm x,r)$ corresponding to $X_N$. We show Figure~\ref{Fig.cap} an example of a spherical cap and the point set inside it.

Additionally, we also can generate index sets of $\bm x\in X_N$ based on \emph{$k$ points nearest neighbor (knn)} $\mathcal C_{K}$ by finding $k$ points $\{\bm x_1,\ldots,\bm x_k\}\in X_N$ in the smallest spherical cap neighbor $C(\bm x,k)\subset\mathbb S^2$, which is
\begin{align*}
C(\bm x,k):=\{\{\bm x_1,\ldots,\bm x_k\}:\exists r>0, \forall\bm x_m\in\{\bm x_1,\ldots,\bm x_k\}, s.t. \lVert\bm x\times\bm x_m\rVert\leq r\}.
\end{align*}
Then we have $k$ number of points $\{\bm x_1,\ldots,\bm x_k\}\in C(\bm x,k)$. Thus, we can obtain index set $\mathcal C_{K}(X_N)$ of points in $C(\bm x,r)$ corresponding to $X_N$.

The neighborhood point sets given by $C(\bm x,r)$ or $C(\bm x,k)$ can be used in fine-tuned denoising for spherical signal processing. 
\subsection{Thresholding techniques}
\label{subsec:thr}
Given a noisy signal $f_\sigma$ with  Gaussian white noise $N(0,\sigma^2)$, we have $f_{\sigma}=f+g$ with $f\in\Pi_{t_{J}}$ and residual $g\notin\Pi_{t_{J}}$. Given the truncated spherical framelet system $\mathcal{F}_{J_0}^{J}(\eta,\mathcal Q)$ as in \eqref{def:sphTruncated}, since $f\in\Pi_{t_{J}}$, then we have $f = \sum_{k=1}^{N_{0}}\mathpzc v_{J_0,k}\varphi_{J_0,k}+\sum_{j=J_0}^J\sum_{k=1}^{N_{j}}\sum_{s=1}^n \mathpzc w_{j,k}^{(s)}\psi_{j,k}^{(s)}$, where coefficients $\mathpzc v_{j,k}:=\langle f,\varphi_{j,k}\rangle$ and  $\mathpzc w_{j,k}^{(s)}:=\langle f,\psi_{j,k}^{(s)}\rangle$.

For $f$, the coefficients need to be normalized first. we compute the energy norm $\lVert \varphi_{j,k}\rVert_{L_2(\mathbb S^2)}$ and $\lVert\psi_{j,k}^{(s)}\rVert_{L_2(\mathbb S^2)}$ for $s=1,\ldots,n$, then the coefficients become
\begin{align*}
\tilde{\mathpzc v}_{j,k}=\frac{\mathpzc v_{j,k}}{\lVert\varphi_{j,k}\rVert_{L_2(\mathbb S^2)}},\quad
\tilde{\mathpzc w}_{j,k}^{(s)}=\frac{\mathpzc w_{j,k}^{(s)}}{\lVert \psi_{j,k}^{(s)}\rVert_{L_2(\mathbb S^2)}},\quad s=1,\ldots,n.
\end{align*}

Second, the coefficients can be processed with various thresholding techniques.  There are four ways for thresholding in this paper:  Global hard  (GH) and global soft (GS) thresholding techniques and local hard (LH) and local soft thresholding techniques. 
We represent $\{\breve{\mathpzc v}_{j,k};\breve{\mathpzc w}_{j,k}^{(s)}:s=1,\ldots,n\}$ as global thresholds and $\{\check{\mathpzc v}_{j,k};\check{\mathpzc w}_{j,k}^{(s)}:s=1,\ldots,n\}$ as local thresholds. 

The global hard threshold (GH)  and the global soft (GS) threshold are given by 
\begin{align}
\label{hrd1}
\breve{\mathpzc w}_{j,k}^{(s)}=
\begin{cases}
\tilde{\mathpzc w}_{j,k}^{(s)},&\lvert\tilde{\mathpzc w}_{j,k}^{(s)}\rvert\geq \tau,\\
0,&\lvert\tilde{\mathpzc w}_{j,k}^{(s)}\rvert< \tau,
\end{cases}
%\label{hrd2}
\breve{\mathpzc v}_{j,k}=
\begin{cases}
\tilde{\mathpzc v}_{j,k},&\lvert\tilde{\mathpzc v}_{j,k}\rvert\geq \tau,\\
0,&\lvert\tilde{\mathpzc v}_{j,k}\rvert< \tau,
\end{cases}
\end{align}
and
\begin{align}
\label{soft1}
\breve{\mathpzc w}_{j,k}^{(s)}=
\begin{cases}
\tilde{\mathpzc w}_{j,k}^{(s)}-\sgn(\tilde{\mathpzc w}_{j,k}^{(s)})\sigma,&\lvert\tilde{\mathpzc w}_{j,k}^{(s)}\rvert\geq \tau,\\
0,&\lvert\tilde{\mathpzc w}_{j,k}^{(s)}\rvert< \tau,
\end{cases}
%\label{soft2}
\breve{\mathpzc v}_{j,k}=
\begin{cases}
\tilde{\mathpzc v}_{j,k}-\sgn(\tilde{\mathpzc v}_{j,k})\sigma,&\lvert\tilde{\mathpzc v}_{j,k}\rvert\geq \tau,\\
0,&\lvert\tilde{\mathpzc v}_{j,k}\rvert< \tau,
\end{cases}
\end{align}
respectively, where $\tau=c\sigma$ is the threshold value and $c$ is a constant. 

The local hard (LH) spherical cap threshold and local soft (LS) spherical cap threshold are introduced below. For each $\bm x_k^j\in X_{N_j}$, let it be the centroid of a spherical cap $C(\bm x_k^j,m)\subset\mathbb S^2$, where $m$ could be radius or number of points. We can find point set $\mathcal X_{M}(\bm x_k^j):=\mathcal X_{M_j^k}\subset C(\bm x_k^j,m)$, where elements $\bm x\in\mathcal X_{M_j^k}$ is selected in $X_{N_{j}}$ by the index of $\mathcal C_{m}(X_{N_j})$, we assume that there are $M_j^k$ points in $\mathcal X_{M_j^k}$. Now we have the coefficients on spherical cap to $\bm x_k^j$, which are $\{\tilde{\mathpzc v}_{j,i}:\bm x_i\in\mathcal X_{M_j^k}\}$ and $\{\tilde{\mathpzc w}_{j,i}^{(1)},\ldots,\tilde{\mathpzc w}_{j,i}^{(n)}:\bm x_i\in\mathcal X_{M_{j+1}^k}\}$, then we averaging each coefficients, which are
\begin{align}
\label{eq36}
\bar{\mathpzc w}_{j,k}^{(s)}&=\frac{1}{M_{j+1}^k}\sum_{i=1}^{M_{j+1}^k} \lvert\tilde{\mathpzc w}_{j,i}^{(s)}\rvert^2,\quad s=1,\ldots,n,\, \bm x_i\in\mathcal X_{M_{j+1}^k}\\
\label{eq37}
\bar{\mathpzc v}_{j,k}&=\frac{1}{M_{j}^k}\sum_{i=1}^{M_{j}^k} \lvert\tilde{\mathpzc v}_{j,i}\rvert^2,\quad\bm x_i\in\mathcal X_{M_{j}^k}.
\end{align}
We define
\begin{align}
\sigma_{j,k}^{(s)}&=\sqrt{\max\{\bar{\mathpzc w}_{j,k}^{(s)}-\sigma^2,0\}},\quad s=1,\ldots,n,\\
\sigma'_{j,k}&=\sqrt{\max\{\bar{\mathpzc v}_{j,k}-\sigma^2,0\}}.
\end{align}
Then we do the threshold value estimation, which is
\begin{align}
\label{eq38}
\tau_{j,k}^{(s)}&=\frac{c\sigma^2}{\sigma_{j,k}^{(s)}},\quad s=1,\ldots,n,\\
\label{eq39}
\tau'_{j,k}&=\frac{c\sigma^2}{\sigma'_{j,k}},
\end{align}
where $c$ is a constant. The local hard (LH) shrinkage is given by
\begin{align}
\label{Lhrd1}
\check{\mathpzc w}_{j,k}^{(s)}&=
\begin{cases}
\tilde{\mathpzc w}_{j,k}^{(s)},&\lvert\tilde{\mathpzc w}_{j,k}^{(s)}\rvert\geq \tau_{j,k}^{(s)},\\
0,&\lvert\tilde{\mathpzc w}_{j,k}^{(s)}\rvert< \tau_{j,k}^{(s)},
\end{cases}\\
\label{Lhrd2}
\check{\mathpzc v}_{j,k}&=
\begin{cases}
\tilde{\mathpzc v}_{j,k},&\lvert\tilde{\mathpzc v}_{j,k}\rvert\geq \tau'_{j,k},\\
0,&\lvert\tilde{\mathpzc v}_{j,k}\rvert< \tau'_{j,k}.
\end{cases}
\end{align}
Similarly,  the shrinkage for local soft (LS) thresholds is given by
\begin{align}
\label{eq40}
\check{\mathpzc w}_{j,k}^{(s)}&=
\begin{cases}
\tilde{\mathpzc w}_{j,k}^{(s)}-\sgn(\tilde{\mathpzc w}_{j,k}^{(s)})\tau_{j,k}^{(s)},&\lvert\tilde{\mathpzc w}_{j,k}^{(s)}\rvert\geq \tau_{j,k}^{(s)},\\
0,&\lvert\tilde{\mathpzc w}_{j,k}^{(s)}\rvert< \tau_{j,k}^{(s)},
\end{cases}\\
\label{eq41}
\check{\mathpzc v}_{j,k}&=
\begin{cases}
\tilde{\mathpzc v}_{j,k}-\sgn(\tilde{\mathpzc v}_{j,k})\tau'_{j,k},&\lvert\tilde{\mathpzc v}_{j,k}\rvert\geq \tau'_{j,k},\\
0,&\lvert\tilde{\mathpzc v}_{j,k}\rvert< \tau'_{j,k}.
\end{cases}
\end{align}

Finally, we need to denormalized the coefficients, which is
\begin{align}
\label{eq42}
\dot{\mathpzc w}_{j,k}^{(s)}&= \mathring{\mathpzc w}_{j,k}^{(s)}\lVert \psi_{j,k}^{(s)}\rVert_{L_2(\mathbb S^2)},\quad s=1,\ldots,n,\\
\label{eq43}
\dot{\mathpzc v}_{j,k}&=\mathring{\mathpzc v}_{j,k}\lVert \varphi_{j,k}\rVert_{L_2(\mathbb S^2)},
\end{align}
where $\{\mathring{\mathpzc v}_{j,k};\mathring{\mathpzc w}_{j,k}^{(1)},\ldots,\mathring{\mathpzc w}_{j,k}^{(n)}\}$ is one of thresholds (global hard, global soft, local hard and local soft) we mentioned above. Thus, we obtain new coefficients $\{\dot{\mathpzc v}_{j,k};\dot{\mathpzc w}_{j,k}^{(1)},\ldots,\dot{\mathpzc w}_{j,k}^{(n)}\}$, which has been denoised by local soft spherical cap threshold. 

Following the procedure of the spherical framelet system, we reconstruct the function with $\{\dot{\mathpzc v}_{j,k};\dot{\mathpzc w}_{j,k}^{(1)},\ldots,\dot{\mathpzc w}_{j,k}^{(n)}\}$, which is
\begin{align}
\label{eq44}
f_{thr}(\bm x)=\sum_{k=1}^{N_j}\dot{\mathpzc v}_{j,k}\varphi_{j,k}+\sum_{k=1}^{N_{j+1}}\sum_{s=1}^{n}\dot{\mathpzc w}_{j,k}^{(s)}\psi_{j,k}^{(s)}.
\end{align}

For $g$, on the point set $\bm x_k^j\in X_{N_j}$, both the hard and soft thresholds are similar to $f$. That is,  the global hard (GH) threshold is given as
\begin{align}
\label{eqg1}
g_{\tau_1}(\bm x_k^j)=
\begin{cases}
g(\bm x_k^j),&\lvert g(\bm x_k^j)\rvert\geq \tau_1,\\
0,&\lvert g(\bm x_k^j)\rvert< \tau_1,
\end{cases}
\end{align}
and the global soft (GS) threshold is given as
\begin{align}
\label{eqg2}
g_{\tau_1}(\bm x_k^j)=
\begin{cases}
g(\bm x_k^j)-\sgn(g(\bm x_k^j))\sigma,&\lvert g(\bm x_k^j)\rvert\geq \tau_1,\\
0,&\lvert g(\bm x_k^j)\rvert< \tau_1,
\end{cases}
\end{align}
where $\tau_1=c_1\sigma$ is the threshold value and $c_1$ is a constant.
The local hard (LH) spherical cap threshold for $\bm x_k^j\in X_{N_j}$ uses $\mathcal X_{M_j^k}\subset C(\bm x_k^j,m)$, which is
\begin{align}
\label{eqg3}
g_{thr}(\bm x_k^j)&=
\begin{cases}
g(\bm x_k^j),&\lvert g(\bm x_k^j)\rvert\geq \tau_{j,k},\\
0,&\lvert g(\bm x_k^j)\rvert< \tau_{j,k},
\end{cases}
\end{align}
where $\tau_{j,k}=\frac{c_1\sigma^2}{\sqrt{\max\{\bar g(\bm x_k^j)-\sigma^2,0\}}}$ and $\bar g(\bm x_k^j)$ with the form
\begin{align}
\bar g(\bm x_k^j)=\frac{1}{M_{j}^k}\sum_{i=1}^{M_{j}^k} \lvert g(\bm x_i^j)\rvert^2,\quad\bm x_i^j\in\mathcal X_{M_{j}^k},
\end{align}
while the local soft (LS) spherical cap threshold for $\bm x_k^j\in X_{N_j}$ uses $\mathcal X_{M_j^k}\subset C(\bm x_k^j,m)$, which is
\begin{align}
\label{eqg4}
g_{thr}(\bm x_k^j)&=
\begin{cases}
g(\bm x_k^j)-\sgn(g(\bm x_k^j))\tau_{j,k},&\lvert g(\bm x_k^j)\rvert\geq \tau_{j,k},\\
0,&\lvert g(\bm x_k^j)\rvert< \tau_{j,k}.
\end{cases}
\end{align}
Hence, we obtain $g_{thr}$ after local soft thresholding.

\subsection{Denoising of Wendland function}
In this section, we give some experiments to illustrate the efficiency of denoising. That is, after projecting noisy function $f_\sigma$ with noise level $\sigma$ on $\Pi_t$ space on spherical $t$-design point sets, we have $f_\sigma=f+g$. We use spherical $t$-design point sets along with its equal weight quadrature rules for framelet decomposition and reconstruction on $\mathbb S^2$. Finally, we apply the thresholds introduced in Section~\ref{subsec:thr} for denoising. Therefore, we have a new function $F_{thr}=f_{thr}+g_{thr}$. 
%To verify the quality of the denoised $F_{thr}$,  we use Signal-to-noise ratio (SNR).  The SNR (with unit dB) is defined as $\text{SNR}:=10\log_{10}(\frac{P_{s}}{P_{n}})$, where $P_{s}$ is the power of a signal and $P_{n}$ is the power of noise, which defined as $P_{\bm y}=(\sum_{i=1}^n y_i^2)^{\frac12}$ for $\bm y:=(y_1,\ldots,y_n)^\top\in\mathbb R^n$. 

Based on the spherical $t$-design point sets (collected by TR-PCG method) in \cite{xiao2023spherical} and the experiments of the best polynomial degrees of projection errors of Wendland functions in Section~\ref{subsec:smth}, we choose spherical $t_j$-design point sets with degree $t_j=16,32,64$ both in SPD and SUD, $t_j=11,24,49$ in SID, and $t_j=12,26,54$ in SHD ($j=1,2,3$). Then we use $X_{N_J}$ to project noisy data $f_\sigma=f_4+G_\sigma$ (generated by normalized Wendland function $f_4$ in \eqref{eq:Phi} and Gaussian white noise $G_{\sigma}$ with noise level $\sigma\lvert f_4\rvert_{\max}$ and $\sigma=[0.05:0.025:0.2]$) on $\Pi_{t_J}$ space with maximum degree $t_J$ ($J=3$) by Algorithm~\ref{alg:projCG} under the setting of maximum iterations $K_{\max}=1000$ and termination tolerance $\varepsilon$ = 2.2204e-16 (floating-point relative accuracy of MATLAB) to obtain $f_\sigma=f+g$, then decompose and reconstruct using a 2-level spherical framelet system $\mathcal{F}_{J_0}^J(\eta,\mathcal Q)$ with four kinds thresholds for denoising: global hard (GH), global soft (GS), local hard (LH) and local soft (LS). 

\begin{figure}[htpb!]
\centering

%\subfigure[3D view]{
\includegraphics[width=0.5\textwidth]{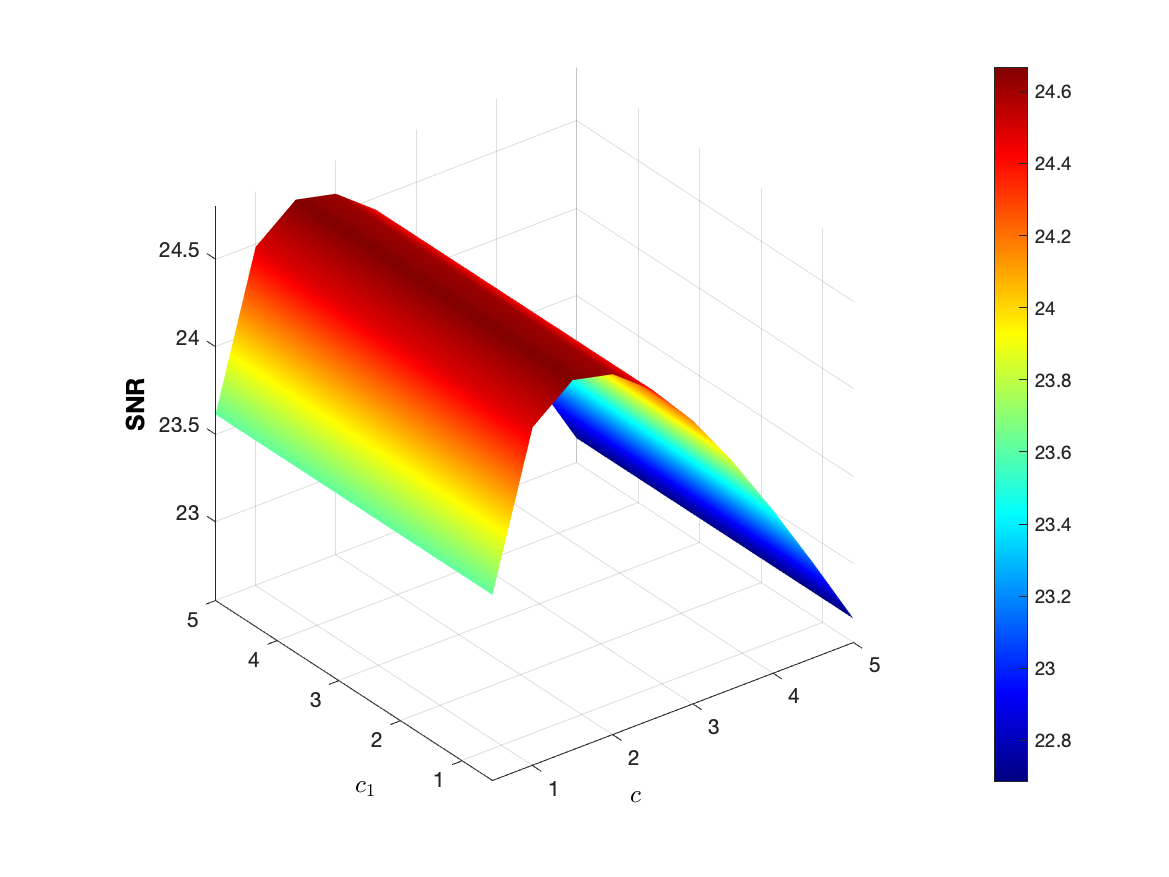}
%}
%\quad
%\subfigure[Top view]{
%\includegraphics[width=0.33\textwidth]{figcc1_wend_top.eps}
%}

\caption{The general view for the behavior under different setting of $c$ and $c_1$ in denoising $f_{0.05}$ by $\eta_3$ on SPD for $t_j=64,32,16$ at cap layer order $i=27$.}
\label{Fig.wend.cc1}
\end{figure}

We use SNR for measuring original signal $f_4$ and the rest of noise $F_{thr}-f_4$, which is $\text{SNR}(f_4,F_{thr}):=10\log_{10}(\frac{\lVert f_4\rVert}{\lVert F_{thr}-f_4\rVert})$ in different filter bank $\eta_k$ for $k=1,2,3$ (see \cite[Section 5.1]{xiao2023spherical}),
 and different point sets $X_{N_j}$ corresponding with degree $t_j$. For finding the suitable constants in Section~\ref{subsec:fmt}, we do a lot experiments by changing the values of $c$ and $c_1$ to see the variation of $\text{SNR}(f_4,F_{thr})$, we give a general view for the behavior of SNR in using $\eta_3$ at $\sigma=0.05$ in Figure~\ref{Fig.wend.cc1}. Then we set $c=2.5,c_1=3$ in GH and LH, $c=1,c_1=3$ in GS and LS. For local thresholds, the spherical cap we take cap layer orders $i=15,22,27$ for filters $\eta_1,\eta_2,\eta_3$, respectively. Notice that all the settings of thresholds may not obtain the best values of $\text{SNR}(f_4,F_{thr})$, but are the mildest values.  We show the results in Table~\ref{table:wen} for different initial point sets of spherical $t$-designs. In tables, we can see: $\text{SNR}_0$ are the $\text{SNR}(f_4,f_{\sigma})$ for different $\sigma$ values; For each thresholding method, the first row to the third row are: $\text{SNR}(f_4,F_{thr})$ values of using $\eta_1,\eta_2,\eta_3$, respectively.

\begin{table}[htpb!]
\centering
\caption{Wendland $f_4$ denoising results}\label{table:wen}
\begin{footnotesize}
\begin{tabular}{cccccccccc}
\hline
$Q_{N_j}$ & Thr & $\sigma$ & 0.05 & 0.075 & 0.1 & 0.125 & 0.15 & 0.175 & 0.2 \\
\hline
\multirow{13}{*}{SPD} & $\text{SNR}_0$ & ~ & \textbf{13.63} & \textbf{10.11} & \textbf{7.61} & \textbf{5.67} & \textbf{4.09} & \textbf{2.75} & \textbf{1.59} \\
\cline{2-10}
~ & \multirow{3}{*}{GH} & $\eta_1$ & 18.51 & 15.72 & 13.90 & 13.08 & 12.50 & 11.99 & 11.51 \\
~ & ~ & $\eta_2$ & 21.50 & 17.64 & 15.11 & 13.89 & 13.23 & 12.74 & 12.28 \\
~ & ~ & $\eta_3$ &  23.25 & 19.13 & 15.75 & 14.31 & 13.55 & 13.09 & 12.67 \\
\cline{3-10}
~ & \multirow{3}{*}{GS} & $\eta_1$ & 19.10 & 16.45 & 14.65 & 13.26 & 12.13 & 11.21 & 10.43 \\
~ & ~ & $\eta_2$ & 20.50 & 17.48 & 15.43 & 13.89 & 12.66 & 11.64 & 10.80 \\
~ & ~ & $\eta_3$ & 21.29 & 18.09 & 15.90 & 14.24 & 12.93 & 11.86 & 10.98 \\
\cline{3-10}
~ & \multirow{3}{*}{LH} & $\eta_1$ & 19.94 & 16.97 & 15.29 & 14.29 & 13.80 & 12.84 & 12.27 \\
~ & ~ & $\eta_2$ & 22.80 & 19.53 & 17.43 & 15.91 & 14.55 & 13.28 & 12.81 \\
~ & ~ & $\eta_3$ & 24.36 & 21.04 & 18.82 & 17.15 & 15.52 & 13.89 & \textbf{13.21} \\
\cline{3-10}
~ & \multirow{3}{*}{LS} & $\eta_1$ & 20.67 & 18.06 & 16.42 & 15.21 & 14.19 & 13.24 & 12.31 \\ 
~ & ~ & $\eta_2$ & 23.11 & 20.05 & 18.03 & 16.47 & 15.18 & 14.02 & 12.88 \\ 
~ & ~ & $\eta_3$ & \textbf{24.48} & \textbf{21.25} & \textbf{19.03} & \textbf{17.30} & \textbf{15.82} & \textbf{14.49} & 13.19 \\ 
\hline\hline
\multirow{13}{*}{SUD} & $\text{SNR}_0$ & ~ & \textbf{13.63} & \textbf{10.11} & \textbf{7.61} & \textbf{5.67} & \textbf{4.09} & \textbf{2.75} & \textbf{1.59} \\
\cline{2-10}
~ & \multirow{3}{*}{GH} & $\eta_1$ & 18.56 & 15.81 & 13.69 & 12.77 & 12.15 & 11.58 & 10.87 \\
~ & ~ & $\eta_2$ & 21.57 & 17.80 & 14.74 & 13.61 & 13.07 & 12.50 & 11.90 \\
~ & ~ & $\eta_3$ & 22.97 & 19.27 & 15.30 & 14.07 & 13.47 & 12.92 & 12.34 \\
\cline{3-10}
~ & \multirow{3}{*}{GS} & $\eta_1$ & 19.05 & 16.29 & 14.47 & 13.07 & 11.91 & 10.96 & 10.17 \\
~ & ~ & $\eta_2$ & 20.44 & 17.40 & 15.30 & 13.71 & 12.45 & 11.43 & 10.58 \\
~ & ~ & $\eta_3$ & 21.15 & 17.92 & 15.73 & 14.07 & 12.75 & 11.68 & 10.81 \\
\cline{3-10}
~ & \multirow{3}{*}{LH} & $\eta_1$ & 19.77 & 16.85 & 15.39 & 14.49 & 13.92 & 13.03 & 12.49 \\
~ & ~ & $\eta_2$ & 22.51 & 19.34 & 17.25 & 15.82 & 14.51 & 13.36 & 12.68 \\
~ & ~ & $\eta_3$ & 23.91 & 20.48 & 18.36 & 16.84 & 15.36 & 13.84 & \textbf{13.04} \\
\cline{3-10}
~ & \multirow{3}{*}{LS} & $\eta_1$ & 20.70 & 18.02 & 16.39 & 15.23 & 14.22 & 13.20 & 12.18 \\
~ & ~ & $\eta_2$ & 22.98 & 19.98 & 17.95 & 16.38 & 15.09 & 13.90 & 12.70 \\
~ & ~ & $\eta_3$ & \textbf{24.15} & \textbf{20.90} & \textbf{18.73} & \textbf{17.07} & \textbf{15.63} & \textbf{14.29} & 12.97 \\
\hline\hline
\multirow{13}{*}{SID} & $\text{SNR}_0$ & ~ & \textbf{13.63} & \textbf{10.11} & \textbf{7.61} & \textbf{5.67} & \textbf{4.09} & \textbf{2.75} & \textbf{1.59} \\
\cline{2-10}
~ & \multirow{3}{*}{GH} & $\eta_1$ & 18.61 & 14.60 & 12.60 & 11.26 & 9.35 & 7.69 & 7.10 \\
~ & ~ & $\eta_2$ & 22.22 & 18.39 & 15.40 & 12.96 & 10.18 & 8.32 & 7.69 \\
~ & ~ & $\eta_3$ & 23.51 & 20.02 & 16.94 & 14.18 & 10.92 & 8.69 & 7.96 \\
\cline{3-10}
~ & \multirow{3}{*}{GS} & $\eta_1$ & 17.07 & 14.04 & 12.04 & 10.57 & 9.41 & 8.43 & 7.58 \\
~ & ~ & $\eta_2$ & 18.62 & 15.37 & 13.12 & 11.42 & 10.06 & 8.94 & 7.98 \\
~ & ~ & $\eta_3$ & 19.25 & 15.94 & 13.62 & 11.84 & 10.42 & 9.24 & 8.26 \\
\cline{3-10}
~ & \multirow{3}{*}{LH} & $\eta_1$ & 19.92 & 16.44 & 14.03 & 12.65 & 11.61 & 10.88 & 10.25 \\
~ & ~ & $\eta_2$ & 23.44 & 19.98 & 17.43 & 15.78 & 14.30 & 13.10 & 11.94 \\
~ & ~ & $\eta_3$ & \textbf{24.66} & \textbf{21.21} & \textbf{18.82} & \textbf{17.00} & \textbf{15.39} & \textbf{14.02} & \textbf{12.79} \\
\cline{3-10}
~ & \multirow{3}{*}{LS} & $\eta_1$ & 20.13 & 16.86 & 14.70 & 13.15 & 11.91 & 10.82 & 9.77 \\ 
~ & ~ & $\eta_2$ & 23.34 & 19.90 & 17.44 & 15.51 & 13.90 & 12.43 & 11.05 \\ 
~ & ~ & $\eta_3$ & 24.54 & 21.03 & 18.51 & 16.46 & 14.68 & 13.09 & 11.59 \\ 
\hline\hline
\multirow{13}{*}{SHD} & $\text{SNR}_0$ & ~ & \textbf{13.58} & \textbf{10.05} & \textbf{7.55} & \textbf{5.62} & \textbf{4.03} & \textbf{2.69} & \textbf{1.53} \\
\cline{2-10}
~ & \multirow{3}{*}{GH} & $\eta_1$ & 17.83 & 14.74 & 12.29 & 10.17 & 8.97 & 8.14 & 7.52 \\
~ & ~ & $\eta_2$ & 21.19 & 17.42 & 14.05 & 11.19 & 9.70 & 8.81 & 8.10 \\
~ & ~ & $\eta_3$ & 22.70 & 18.63 & 15.29 & 11.95 & 10.09 & 9.16 & 8.41 \\
\cline{3-10}
~ & \multirow{3}{*}{GS} & $\eta_1$ & 17.02 & 14.04 & 12.03 & 10.52 & 9.32 & 8.33 & 7.51 \\
~ & ~ & $\eta_2$ & 18.36 & 15.11 & 12.88 & 11.21 & 9.87 & 8.79 & 7.91 \\
~ & ~ & $\eta_3$ & 18.93 & 15.61 & 13.31 & 11.57 & 10.18 & 9.04 & 8.11 \\
\cline{3-10}
~ & \multirow{3}{*}{LH} & $\eta_1$ & 19.54 & 15.86 & 14.02 & 12.60 & 11.40 & 10.32 & 9.34 \\
~ & ~ & $\eta_2$ & 22.14 & 18.56 & 16.06 & 14.40 & 13.09 & 11.77 & 10.55 \\
~ & ~ & $\eta_3$ & \textbf{23.23} & 19.69 & 17.16 & 15.25 & 13.85 & \textbf{12.55} & \textbf{11.25} \\
\cline{3-10}
~ & \multirow{3}{*}{LS} & $\eta_1$ & 19.91 & 16.68 & 14.64 & 13.07 & 11.76 & 10.60 & 9.50 \\ 
~ & ~ & $\eta_2$ & 22.21 & 18.80 & 16.42 & 14.57 & 13.00 & 11.59 & 10.28 \\ 
~ & ~ & $\eta_3$ & 23.21 & \textbf{19.70} & \textbf{17.22} & \textbf{15.26} & \textbf{13.59} & 12.07 & 10.64 \\ 
\hline
\end{tabular}
\end{footnotesize}
\end{table}

\begin{figure}[htpb!]
\centering

\subfigure[$f_4$]{
\includegraphics[width=0.25\textwidth]{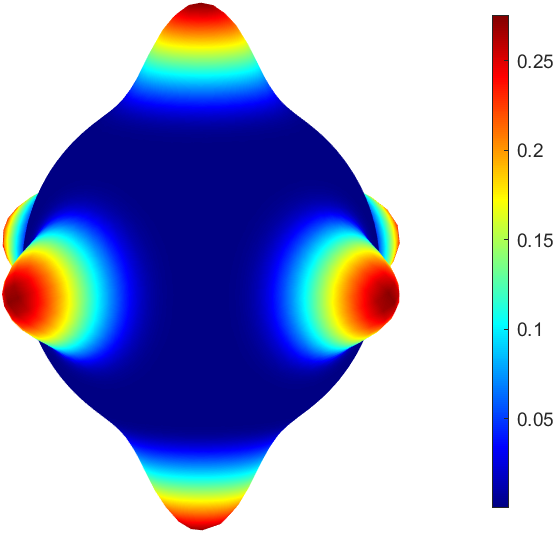}
}
\quad
\subfigure[$f_{0.05}$]{
\includegraphics[width=0.25\textwidth]{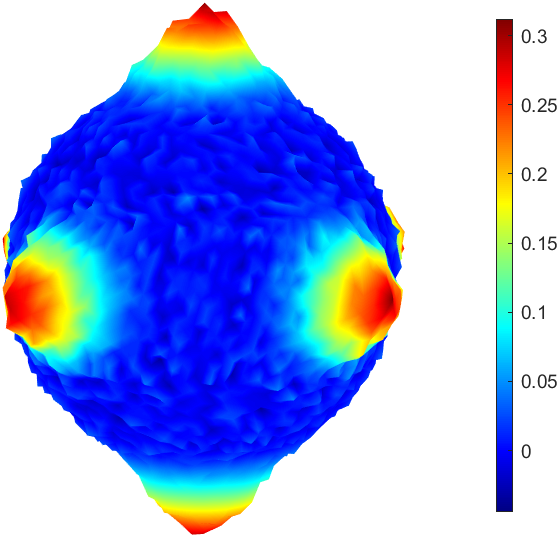}
}
\quad
\subfigure[$f$]{
\includegraphics[width=0.25\textwidth]{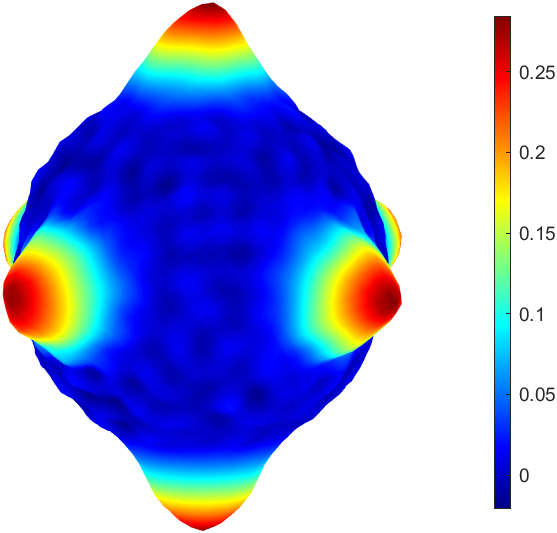}
}
\quad
\subfigure[$g$]{
\includegraphics[width=0.25\textwidth]{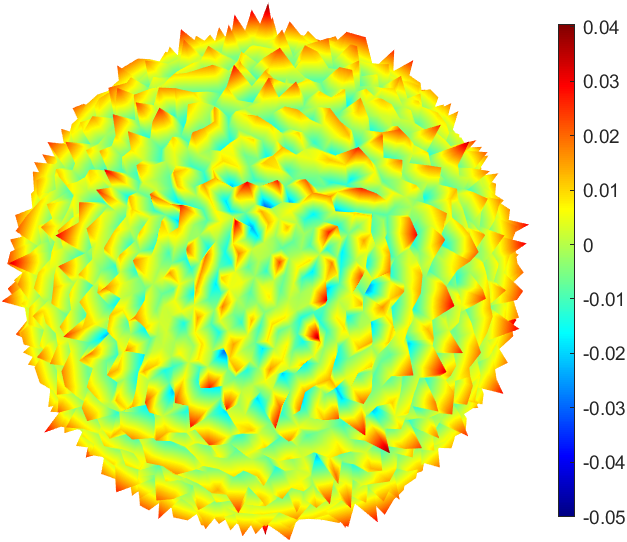}
}
\quad
\subfigure[$F_{thr}$]{
\includegraphics[width=0.25\textwidth]{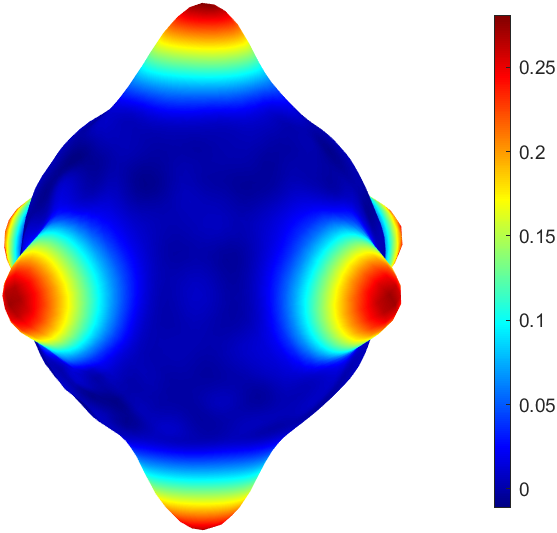}
}
\quad
\subfigure[$f_4-F_{thr}$]{
\includegraphics[width=0.25\textwidth]{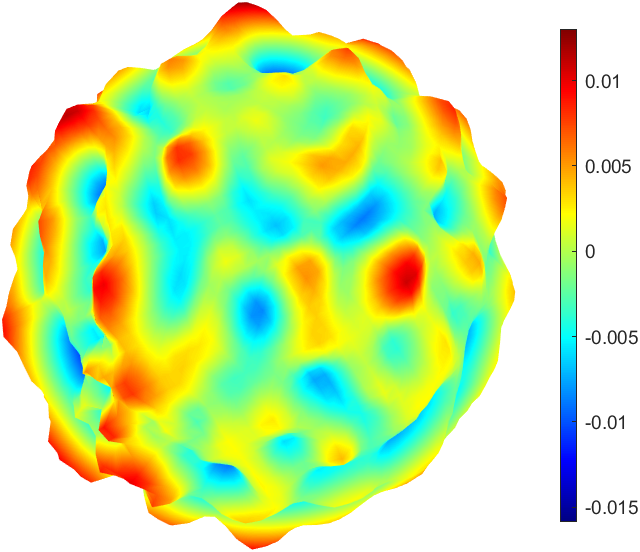}
}

\caption{The behavior of denoising Wendland $f_{0.05}$ by $\eta_3$ with LS on SPD with $t_j=16,32,64$ }
\label{Fig.sigwend}
\end{figure}

As we can see from the tables, local threshold methods (especially local soft (LS)) are more beneficial than global thresholds, which are reasonable in construction. The behavior of different filters is $\eta_3>\eta_2>\eta_1$, which means the more high pass filters, the better behavior in denoising. One interesting discovery is that the initial icosahedron point sets of spherical $t$-designs (SID) starting from $t=49$ (the smallest degree as finest scale) are adapted with the local hard threshold (LH) with good performance for noise level with $\sigma=0.05$ to $0.1$.  Moreover, the initial spiral point sets of spherical $t$-designs (SPD) are reliable in more general applications since we can generate it for what we require with various relations of $t$ and $N$. Figur~\ref{Fig.sigwend} shows the behavior of denoised and reconstructed functions.

%%%%%%%%%%%%%%%%%%%%%%%%%%%%%%%%%%%%%%
%%%%%%%%%%%%%%%%%%%%%%%%%%%%%%%%%%%%%%
%%%%%%%%%%%%%%%%%%%%%%%%%%%%%%%%%%%%%%
\section{Conclusions}
\label{sec:conclusions}
In this paper, we present the line search with the restart conjugate gradient method (LS-RCG) for obtaining numerical spherical $t$-design point sets. We also compare the different approaches (LS-RCG and TR-PCG, full Hessian $\mathcal H$ and approximation Hessian $\mathcal H_2$) from the results of their numerical spherical $t$-designs. In general, the LS-RCG consumes a lot  in the line search iterations, whereas the TR-PCG costs fewer iterations in finding the trust regions. However, there are no significant differences in the accuracy of full Hessian and approximation Hessian. Then, we study the approximation of spherical functions (Wendland functions with certain smoothness) and non-smooth functions (spherical $1_{\mathcal D}$ function) by projecting them with orthogonal polynomial using spherical $t$-designs point sets on $\mathbb S^2$. The results show that choosing different weights and different types of spherical $t$-designs will affect the error accuracy of projects. Moreover, we introduce the spherical signal processing based on the best projection degree of polynomials in Wendland functions. We use spherical framelets in \cite{xiao2023spherical} from spherical designs under four different types of thresholding techniques which the fine-tuned spherical caps may demand. We show that all the point sets perform well under the spherical framelets in the denoising procedure, indicating the spherical designs have great potential in spherical data processing.

\begin{acknowledgements}
This work was supported in part by the Research Grants Council of the Hong Kong Special Administrative Region, China, under
Project CityU 11309122 and Project CityU 11302023.
\end{acknowledgements}

\section*{Conflict of interest}
The authors declare that they have no conflict of interest.

\bibliographystyle{spmpsci}
%\bibliography{references}
\end{document}